\documentclass[final,12pt]{elsarticle}

\usepackage{amssymb}
\usepackage{amsmath}
\usepackage{amsthm,}
\usepackage{mathtools}

\usepackage{comment}

\usepackage{booktabs}

\usepackage[dvipsnames]{xcolor}

\usepackage{algorithm}
\usepackage{algpseudocode}

\usepackage{tikz}
\usepackage{caption}
\usepackage{graphicx}
\usepackage{subcaption} 

\usepackage{cases}

\usepackage{stmaryrd}

\begin{document}

\newcommand{\vertiii}[1]{{\left\vert\kern-0.25ex\left\vert\kern-0.25ex\left\vert #1 
    \right\vert\kern-0.25ex\right\vert\kern-0.25ex\right\vert}}
\newcommand{\af}[1]{\textcolor{red}{#1}}
\newcommand{\dm}[1]{\textcolor{blue}{#1}}

\newcommand{\gt}{\boldsymbol{g}_T}
\newcommand{\disp}{\mathbf{u}}
\newcommand{\trac}{\mathbf{t}}
\newcommand{\jump}[1]{\llbracket #1 \rrbracket}
\newcommand{\card}[1]{\vert #1 \vert}
\newcommand{\idx}[2]{\mathcal{I}_{#1}^{(#2)}}

\newcommand{\filledsquare}[2][black]{%
  \raisebox{1.5#2}{%
    \fcolorbox{black}{#1}{\color{#1}\rule{#2}{#2}}%
  }%
}

\newcommand{\filledtriangle}[2][black]{%
  \raisebox{-0.4#2}{%
    \tikz\draw[black, fill=#1, line width=0.3pt] 
      (0,0) -- (#2,0) -- (#2/2,0.87*#2) -- cycle;
  }%
}

\definecolor{stick}{rgb}{0.267,0.0487,0.329}
\definecolor{slip}{rgb}{0.208,0.718,0.472}

\theoremstyle{definition}
\newtheorem{lemma}{Lemma}[section]
\newtheorem{remark}{Remark}[section]
\newtheorem{proposition}{Proposition}[section]
\bibliographystyle{elsarticle-num}

\begin{frontmatter}


\title{A Stabilized Mortar Method for Discontinuities in Geological Media with Non-Conforming Grids}


\affiliation[UniPD]{organization={Dept. of Civil, Environmental and Architectural Engineering, University of Padova},
            state={Padova},
            country={Italy}}

\author[UniPD]{Daniele Moretto\corref{coraut}} 
\ead{daniele.moretto.3@phd.unipd.it}
\author[UniPD]{Andrea Franceschini}
\author[UniPD]{Massimiliano Ferronato}

\cortext[coraut]{Corresponding author}

\begin{abstract}
Accurate numerical simulation of fault and fracture mechanics is critical for the performance and safety assessment of many subsurface systems. 
The discretized representation of discontinuity surfaces and the robust simulation of their frictional contact behavior still represent major challenges.
In this work, we use the mortar method to enforce the contact constraints and allow for non-conformity around the discontinuity surface, with a set of Lagrange multipliers playing the role of interface tractions. The formulation combines piecewise linear displacements with piecewise constant multipliers defined on one side of the
fault interface (the non-mortar side). This choice for the Lagrange multipliers has a number of advantages from practical and computational viewpoints, but violates the inf-sup stability constraint. In order to stabilize the proposed formulation, we develop a traction-jump stabilization term to be added to the constraint equations. We use a macro-element analysis to derive an algorithmic strategy that automatically evaluates the proper scaling of the stabilization, without requiring any additional user-selected parameter. Numerical experiments demonstrate that the proposed formulation not only restores the inf-sup stability condition, but also recovers stable traction profiles in the presence of finer non-mortar sides, where other inf-sup-stable formulations fail. 
The proposed method is finally used to simulate non-linear contact conditions in non-conforming corner-point grids typically used in industrial geological applications.

\end{abstract}



\begin{keyword}
Contact mechanics \sep Lagrange multipliers \sep Mortar method \sep Stabilization \sep Corner-point grid




\end{keyword}
\end{frontmatter}



\section{Introduction}
Numerical modeling of the subsurface geomechanical response plays a crucial role in assessing the safety of many engineering applications, such as geothermal energy production, underground gas storage, and carbon sequestration \cite{urbancic1993microseismicity,ferronato2010geomechanical,dockrill2010structural,wang2016safety,buijze2019review,kumar2023comprehensive}. To obtain accurate predictions,
geomechanical models must properly represent faults and fractures within the domain, with the aid of proper contact mechanics algorithms,
%
as these structures can strongly impact the overall mechanical behavior of the geological media under examination. Despite the extensive literature on contact mechanics and related challenges \cite{johnson1987contact,wriggers2003computational}, the development of robust strategies for fault and fracture modeling in geological problems remains an open and highly debated issue.

There are numerous approaches to simulate geological faults \cite{berre2019flow}, and a comprehensive review is beyond the scope of this paper. If we restrict ourselves to a coarse classification, a major distinguishing feature concerns how faults and fractures are represented in the model. In this work, we focus on the class of Discrete Fracture Models (DFM) \cite{puso2004mortar,franceschini2016novel,garipov2016discrete}, which explicitly discretize a mechanical discontinuity as a lower-dimensional grid conforming to the surrounding background grid. Alternative approaches include Embedded Discrete Fracture Models (EDFM) \cite{oliver2006comparative,borja2008assumed,cusini2021simulation}, where the fault is treated as a separate geometrical object embedded in the background grid, and implicit methods, where an explicit representation is avoided. Examples in this class of techniques include the Phase Field approach \cite{chen2002phase,cajuhi2018phase}, where the discontinuity is described through a continuous and smooth additional variable, and Peridynamics \cite{ye2026peridynamics}, which is based on integral equations instead of classical differential equations.

Within the class of DFM approaches, a further classification is based on the method used to enforce contact constraints. Penalty methods are the simplest strategy; however, the constraint is enforced only approximately, and often this can lead to severely ill-conditioned problems. In contrast, Lagrange multiplier methods (LMM) provide a more accurate enforcement of the constraints \cite{papadopoulos1998lagrange, franceschini2022scalable}, at the cost of introducing an additional set of degrees of freedom on the fault. Augmented Lagrange multiplier techniques \cite{rockafellar1974augmented,curnier1988generalized,frigo2025robust} and the Nitsche method \cite{sanders2012embedded,chouly2013nitsche} represent other notable alternatives that combine features of penalty methods and LMM.

In this work, we adopt a DFM approach combined with LMM to model fault and fracture mechanics in geological media. To mitigate the aforementioned geometrical challenges, we employ the mortar method \cite{bernardi1993domain,belgacem1999mortar,popp2012mortar}, which relaxes the conformity requirement and allows for non-matching grids on the two sides of the fault, a situation that commonly arises when using industry-standard grid formats, such as corner-point grids \cite{ponting1989corner}. We consider the original single-pass mortar formulation \cite{puso2004mortar,tur2009mortar}, where the multiplier space is defined on one side of the interface (\textit{non-mortar}), while the other side is referred to as \textit{mortar}.
The resulting weak formulation presents several challenges. First, it leads to a mixed problem that requires careful choice of the functional spaces to satisfy the inf-sup stability condition \cite{boffi2013mixed}. Moreover, the stability of the Lagrange multipliers is affected by the relative refinement of the mesh on the two sides of the interface \cite{pitkaranta1980local,moretto2025stabilized}. In particular, a finer non-mortar side enriches the multiplier space and may produce an unacceptable oscillatory behavior even when adopting inf-sup stable formulations, such as dual nodal multipliers \cite{wohlmuth2002comparison,popp2010dual} or bubble enrichments \cite{hauret2007discontinuous,frigo2025robust}. 
Finally, classical mortar implementation issues concerning the enforcement of boundary conditions and the treatment of interface cross-points should also be considered \cite{puso2003mesh,zhou2020three}.

To address these issues, this work introduces a novel stabilized mortar formulation combining  piecewise linear displacement fields with piecewise-constant Lagrange multipliers and a traction-jump stabilization term added to the constraint equation. Building upon macro-element theory \cite{elman1996iterative,franceschini2020algebraically,camargo2021macroelement}, we propose an algorithm to properly scale the stabilization contribution without requiring any additional user-selected parameters. 
The manuscript is organized as follows. Section 2 introduces the model problem, while Section 3 presents its discrete mortar formulation. In Section 4, we introduce the stabilization strategy and numerically verify the satisfaction of the inf-sup condition. Section 5 reports several numerical examples to validate the model and assess the effectiveness of the proposed stabilization, including an application to a geological medium discretized with a non-conforming corner-point grid. Finally, Section 6 concludes the paper with some closing remarks.

\section{Problem statement}
\begin{figure}
    \centering
    \begin{subfigure}[b]{0.5\textwidth}
        \centering
        \includegraphics[width=\textwidth]{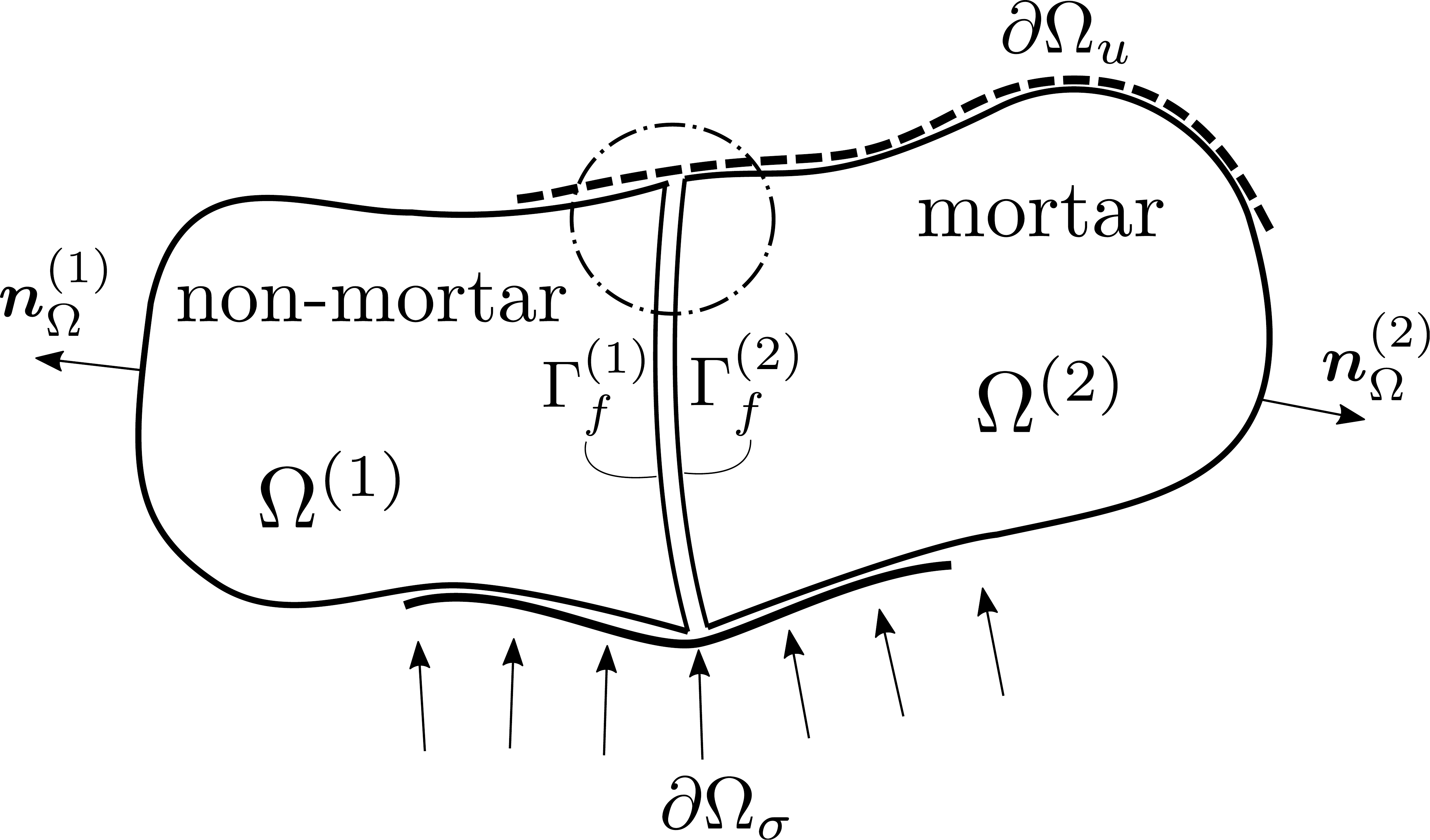}
        \caption{}
        \label{fig:conceptual_scheme}
    \end{subfigure}
    \hspace{0.1\textwidth}
    \begin{subfigure}[b]{0.25\textwidth}
        \centering
        \includegraphics[width=\textwidth]{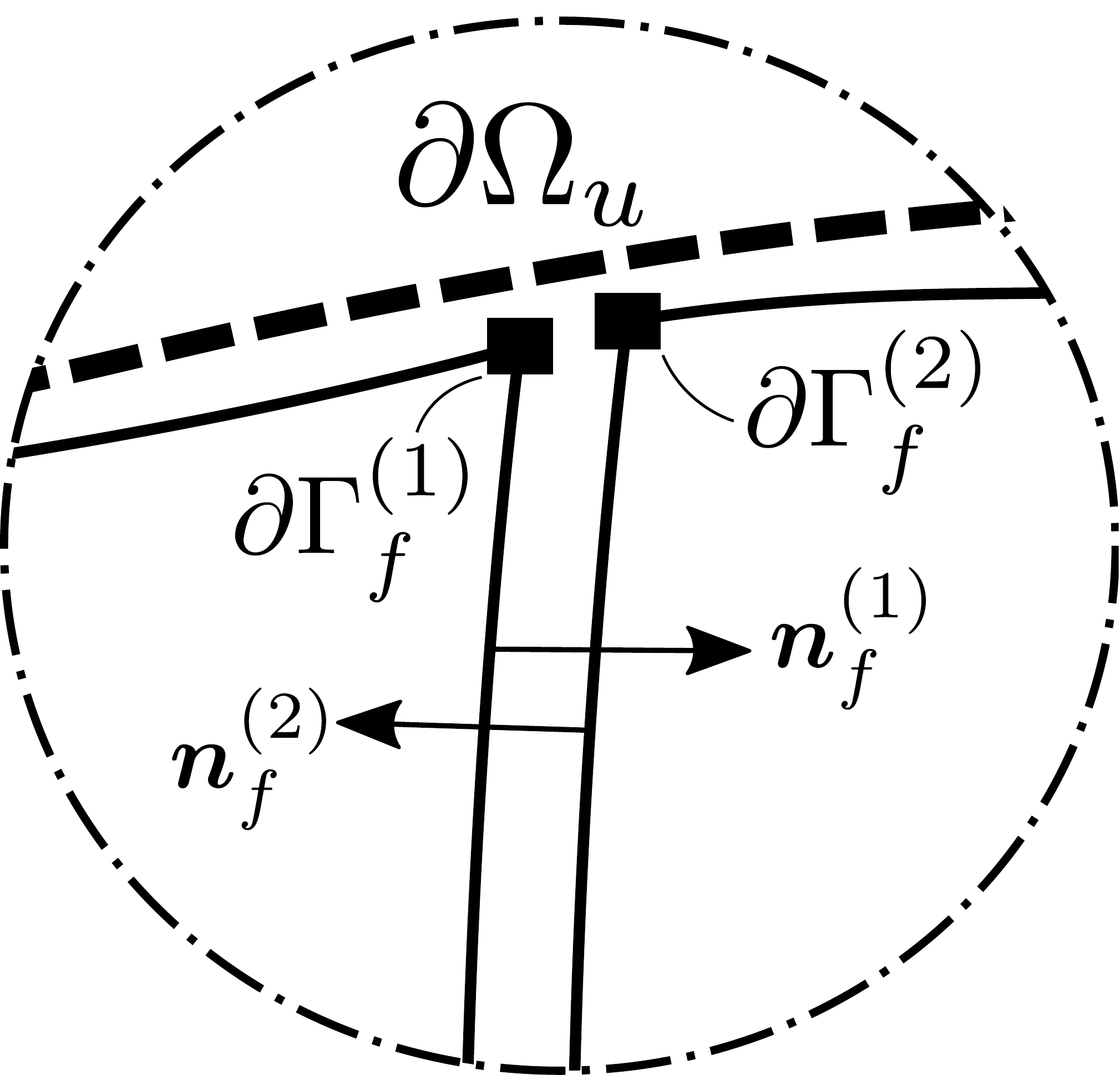}
        \caption{}
        \label{fig:conceptual_scheme_zoom}
    \end{subfigure}
    \caption{Conceptual scheme of the model problem (a) and a detail on the interface boundary (b).}
    \label{fig:conceptual_scheme_main}
\end{figure}

Consider two closed elastic domains $\overline{\Omega}^{(1)}$ and $\overline{\Omega}^{(2)}$ $\subset \mathbb{R}^3$ (see Figure \ref{fig:conceptual_scheme_main}a for reference).
For each domain, we denote by $\Omega^{(i)}$ the interior, $\partial\Omega^{(i)}$ the boundary, and $\boldsymbol n_\Omega^{(i)}$ the outward unit normal. We denote by $\overline{\Gamma}_f^{(i)}$ the closed portion of $\partial\Omega^{(i)}$ that is subject to the contact constraint, and by $\partial\Gamma_f^{(i)}$ its lower-dimensional boundary.
A fault will be represented by the interaction between the non-mortar side $\overline{\Gamma}_f^{(1)}$ and the mortar side $\overline{\Gamma}_f^{(2)}$, which may be geometrically non-conforming, so that a single lower-dimensional inclusion cannot be identified.
Let us also denote by $\boldsymbol n_{f}^{(i)}$ the outward unit normal of $\overline{\Gamma}_f^{(i)}$.
The quantities defined on $\overline{\Gamma}_f^{(i)}$ and oriented along its normal and tangential directions are indicated with subscripts $N$ and $T$, so that  $(\bullet)_N=\boldsymbol{n}_f^{(i)}\cdot(\bullet)$ and $(\bullet)_T = (\mathbf I - \boldsymbol{n}_f^{(i)} \otimes \boldsymbol{n}_f^{(i)})\cdot (\bullet)$ with $\mathbf I$ the identity in $\mathbb{R}^3$ and $\otimes$ the dyadic product.
In the following, contact-related quantities, such as tractions and gaps, will be computed only on the non-mortar side $\Gamma_f^{(1)}$.

Let $\partial\Omega_u$ and $\partial\Omega_{\sigma}$ be the Dirichlet and Neumann subsets, respectively, such that $\partial\Omega_u \cap \partial\Omega_{\sigma} = \emptyset$. We also require that $\Gamma_f^{(1)} \cap (\partial\Omega_u \cup \partial\Omega_{\sigma}) = \emptyset$, but we admit that $\partial \Gamma_f^{(1)} \cap (\partial\Omega_u \cup \partial\Omega_{\sigma}) \neq \emptyset$, that is, the Dirichlet and Neumann boundary can intersect the fault boundary (see Figure \ref{fig:conceptual_scheme_main}b).

In this model, we are interested in simulating the activation of a geological fault or a fracture with no regard for seismicity.
For this reason, we assume infinitesimal strains and negligible inertial effects, 
with the analysis spanning the time interval $[0,t_{\max}]$, for some final instant $t_{\max}>0$.
With these assumptions, the strong form of the frictional contact initial boundary value problem reads: for each $i = \{1,2\}$, given the body forces $\boldsymbol{b}^{(i)}:\Omega^{(i)}\times]0,t_{\max}]\rightarrow\mathbb{R}^3$ and the traction field $\bar{\boldsymbol{t}}:\partial\Omega_\sigma\times]0,t_{\max}]\rightarrow \mathbb{R}^3$, find the displacement field $\boldsymbol{u}^{(i)}: \Omega^{(i)}\times]0,t_{\max}] \rightarrow \mathbb{R}^3$ such that:
\begin{subequations}
    \begin{align}
    &-\nabla \cdot \boldsymbol{\sigma}(\boldsymbol{u}^{(i)})= \boldsymbol{b}^{(i)}, \ &\text{in} \; \Omega^{(i)}\times]0,t_{\max}], \\
    & \boldsymbol u^{(i)} = \boldsymbol 0, \ \ \ &\text{on} \; (\partial\Omega_u \cap \partial\Omega^{(i)})\times]0,t_{\max}],\\
    & \boldsymbol{\sigma}(\boldsymbol{u}^{(i)})\cdot \boldsymbol n_\Omega^{(i)} = \bar{\boldsymbol t}, \ \ \ &\text{on} \; (\partial\Omega_{\sigma} \cap \partial\Omega^{(i)})\times]0,t_{\max}],\\
    & \boldsymbol u^{(i)} = \boldsymbol 0, \ \ \ &\text{in} \; \overline{\Omega}^{(i)} \text{at } t=0,
\end{align}
\label{eq:strongBVP}\null
\end{subequations}
subject to the following contact condition over $\Gamma_f^{(1)}\times]0,t_{\max}]$:
\begin{align}
    &t_N = \boldsymbol{t}\cdot \boldsymbol{n}_f^{(1)} \leq 0, \   g_N = \jump{\boldsymbol u} \cdot \boldsymbol{n}_f^{(1)} \geq 0, \  t_N g_N = 0, \  &\text{(normal contact conditions)}, \nonumber\\
    &\|\boldsymbol{t}_T \|_2 \leq \tau_{lim}(t_N),  \ \dot{\boldsymbol{g}}_T \cdot \boldsymbol{t}_T - \tau_{lim}(t_N)\| \dot{\boldsymbol{g}}_T \|_2 = 0, \ &\text{(quasi-static friction law)}. \nonumber
\end{align}
Homogeneous essential and initial conditions have been considered for convenience.
In \eqref{eq:strongBVP}, $\boldsymbol{\sigma}(\boldsymbol u^{(i)}) = \mathbb{D}: \nabla^s \boldsymbol{u}^{(i)}$ is the stress tensor, with $\mathbb{D}$ the fourth-order elastic constitutive tensor and $\nabla^s$ the symmetric gradient operator, $\boldsymbol t = \boldsymbol \sigma(\boldsymbol u^{(1)}) \cdot \boldsymbol n_f^{(1)}$ is the traction vector acting on $\overline\Gamma_{f}^{(1)}$, and $\boldsymbol g = \jump{\boldsymbol u}$ is the classical gap function, computed with the jump operator $\jump{\cdot}$ defined for a general function $v:\Omega^{(i)}\rightarrow\mathbb{R}$ as
\begin{align}
    \jump{v} = v^{(1)}\vert_{\Gamma_f^{(1)}} - \Pi \big( v^{(2)}\vert_{\Gamma_f^{(2)}} \big),
    \label{eq:interf_jump}
\end{align}
where $\Pi$ is a suitable map that projects $v$ from the mortar to the non-mortar side.
Among the several friction laws available, we adopt the Coulomb criterion, which is widely used in large-scale geological applications mainly because of its simplicity (it requires only two material parameters) and effectiveness in capturing the essential physics of frictional behavior of rocks \cite{jaeger2007fundamentals}. The Coulomb criterion provides $\tau_{lim} = c - t_N \tan(\phi)$, where $c$ is the cohesion and $\phi$ is the friction angle.
We will consider a static version of the frictional law where the velocity $\dot{\boldsymbol{g}}_T$ is replaced with the displacement increment \cite{wohlmuth2011variationally}. Hence, when adopting an implicit time marching scheme, at time $t_n$, $\dot{\boldsymbol{g}}_T$ is replaced by the backward difference operator $\Delta_n \boldsymbol{g}_T$, defined as
\begin{align}
    \Delta_n\boldsymbol{g}_T = (\boldsymbol{g}_T)_n - (\boldsymbol{g}_T)_{n-1}, 
\end{align}
with $(\cdot)_n$ denoting a quantity at time $t_n$. 
The contact constraints yield a partition of $\Gamma_f^{(1)}$ into three disjoint regions:
\begin{itemize}
    \item a \textit{stick} region ($\Gamma_f^{\text{stick}}$), where $t_N < 0$ and the contact conditions are satisfied;
    \item  a \textit{slip} region ($\Gamma_f^{\text{slip}}$), where $t_N < 0$ 
    and the tangential traction reaches the limiting value prescribed by the Coulomb criterion
    \begin{align}
        \boldsymbol t_T\vert_{\Gamma_f^{\text{slip}}} = \boldsymbol t_T^* = \tau_{lim}(t_N)\frac{\Delta_n\boldsymbol g_T}{\|\Delta_n\boldsymbol g_T\|_2},
    \end{align}
    thus allowing for relative slip between $\Gamma_f^{(1)}$ and $\Gamma_f^{(2)}$;
    \item an \textit{open} region ($\Gamma_f^{\text{open}}$), where $t_N > 0$ and no contact takes place.
\end{itemize}

For further details on the strong formulation of the contact problem, we refer the reader to \cite{popp2010dual,popp2012dual,tur2009mortar,wohlmuth2011variationally,yang2005two}.

\section{Numerical model}

\subsection{Discrete variational formulation} \label{sec:discrete_var_setting}
We present the discrete variational setting of the constrained initial boundary value problem \eqref{eq:strongBVP} according to the mortar framework \cite{belgacem1999mortar,bernardi1993domain,wohlmuth2001discretization}. 
%
We use the notation $(\cdot,\cdot)_{\Omega^{(i)}}$ and $\langle \cdot,\cdot\rangle_{\Gamma_f^{(1)}}$ to denote the $L^2$-inner product in $\Omega^{(i)}$ and $\Gamma_f^{(1)}$, respectively, and $\|\cdot\|_{s;\Omega^{(i)}}$ to indicate the natural norm in $H^s(\Omega^{(i)})$.
Finally, the symbol $\card\cdot$ denotes the cardinality of a given set.

Within each domain, we introduce the continuous space ${\boldsymbol{\mathcal V}}^{(i)} = [H^1_0(\Omega^{(i)})]^3 $ and its trace on the non-mortar side $\boldsymbol{\mathcal W} = [H^{1/2}_{00}(\Gamma_f^{(1)})]^3$, whose dual counterpart is denoted by $\boldsymbol{\mathcal{W}}'$.
We also introduce $\boldsymbol{\mathcal{V}} = \boldsymbol {\mathcal{V}}^{(1)}\times \boldsymbol{\mathcal{V}}^{(2)}$, equipped with the broken norm 
$\|\boldsymbol{v}\|^2_{\boldsymbol{\mathcal{V}}} = \|\boldsymbol{v}^{(1)}\|^2_{1;\Omega^{(1)}} + \|\boldsymbol{v}^{(2)}\|^2_{1;\Omega^{(2)}}$ 
for $\boldsymbol{v} \coloneqq (\boldsymbol{v}^{(1)},\boldsymbol{v}^{(2)}) \in \boldsymbol{\mathcal{V}}$.
We assume that the operator $\Pi$, introduced in \eqref{eq:interf_jump}, is smooth enough so that if $\boldsymbol{v} \in \boldsymbol{\mathcal{V}}$, then $\jump{\boldsymbol v} \in \boldsymbol{\mathcal{W}}$. The contact at the interface between $\Omega^{(1)}$ and $\Omega^{(2)}$ is enforced in a weak sense by introducing a set of Lagrange multipliers on the non-mortar side $\Gamma_f^{(1)}$ lying in the space:

\begin{align}
    \boldsymbol{\mathcal{M}}(t_N) = \{\boldsymbol \mu \in \boldsymbol{\mathcal{W}}': \langle \boldsymbol \mu, \boldsymbol v \rangle_{\Gamma_f^{(1)}} \leq \langle \tau_{lim}, \|\boldsymbol v_T\|_2\rangle_{\Gamma_f^{(1)}}, \, \boldsymbol v \in \boldsymbol{\mathcal{W}}, \, \text{with} \; v_N \leq 0 \},
\end{align}
endowed with the norm $\|\boldsymbol \mu\|_{\boldsymbol{\mathcal{M}}} = \|\boldsymbol \mu\|_{-\frac{1}{2};\Gamma_f^{(1)}}$.

Consider now a partition $\mathcal{T}^{(i)}$ of $\Omega^{(i)}$ such that $\Omega^{(i)} =  \bigcup_{K \in \mathcal{T}^{(i)}}  K$ and denote by $h^{(i)}$ its characteristic size. 
In the following numerical experiments, $\mathcal{T}^{(i)}$ will consist exclusively of hexahedrons, but the proposed formulation can also be extended to tetrahedrons without modifications.
We will use the ratio $r_h = h^{(1)}/h^{(2)}$ to measure the relative refinement of the mortar and non-mortar grids.
Finally, $\mathcal{T} = \mathcal{T}^{(1)}\bigcup\mathcal{T}^{(2)}$ denotes the union of the partitions of $\Omega^{(1)}$ and $\Omega^{(2)}$.

Let $\mathcal{F}$ be the set of faces belonging to $\Gamma_f^{(1)}$, such that $\Gamma^{(1)}_f = \bigcup_{F \in \mathcal{F}} F$, and let $\mathcal{E} = \bigcup_{E_{KL} \in \mathcal{E}} E_{KL}$ be the set of internal edges in $\Gamma_f^{(1)}$ such that $E_{KL} =F_K \cap F_L$. 
Let also $\mathcal{N}_\Omega$ be the joint set of nodes in $\Omega^{(1)}$ and $\Omega^{(2)}$.
Moreover, we denote by $\mathcal{I}_u$ the corresponding set of displacement indices, with $\card{\mathcal{I}_u} = 3\card{\mathcal{N}_\Omega}$ and $|\mathcal{I}_u|=\sum_i|\mathcal{I}_u^{(i)}|$. Finally, let $\mathcal{I}_t$ be the set of traction indices in $\Gamma_f^{(1)}$.


Given the finite-dimensional subspaces $\boldsymbol{\mathcal{V}}^{(i)}_h \subset \boldsymbol{\mathcal{V}}^{(i)}$ and $\boldsymbol{\mathcal{M}}_h(t_{N,h}) \subset \boldsymbol{\mathcal{M}}(t_N)$, 
with $\boldsymbol{\mathcal{V}}_h=\boldsymbol{\mathcal{V}}^{(1)}_h\times\boldsymbol{\mathcal{V}}^{(2)}_h$,
the finite element solution to the initial boundary value problem \eqref{eq:strongBVP} reads: for any instant $t\in]0,t_{\max}]$, find $\{ \boldsymbol u_h, \boldsymbol t_h\} \in \boldsymbol{\mathcal{V}}_h \times \boldsymbol{\mathcal{M}}_h(t_{N,h})$ such that:
\begin{subequations}
\begin{align}
    \mathcal{R}_u &= \sum_i(\nabla^s\boldsymbol{\eta}^{(i)},\mathbb{D} : \nabla^s\boldsymbol{u}^{(i)}_h)_{\Omega^{(i)}} + \langle \jump{\boldsymbol \eta}, \boldsymbol t_h\rangle_{\Gamma_f^{(1)}} - \sum_i(\boldsymbol{\eta}^{(i)},\boldsymbol{b}^{(i)})_{\Omega^{(i)}} - \nonumber \\
    & \quad \; \sum_i\langle \boldsymbol{\eta}^{(i)}, \overline{\boldsymbol t} \rangle_{\partial\Omega_\sigma \cap \partial\Omega^{(i)}} = 0, 
    \label{eq:displacement_equation} \\
    \mathcal{R}_t &= \langle t_{N,h} - \mu_N, g_{N,h} \rangle_{\Gamma_f^{(1)}} + \langle \boldsymbol t_{T,h} - \boldsymbol\mu_T, \boldsymbol g_{T,h} \rangle_{\Gamma_f^{(1)}}  \geq 0,
    \label{eq:constraint_inequality}
\end{align}
\label{eq:Discrete_variational_contact}
\end{subequations}
for all $\{\boldsymbol \eta,\boldsymbol \mu\} \in \boldsymbol{\mathcal{V}}_h \times \boldsymbol{\mathcal{M}}_h(t_{N,h})$.

\subsection{Solution strategy} \label{sec:lin_mortar}
In order to resolve the inequality constraint \eqref{eq:constraint_inequality}, we employ an active-set strategy, which consists of iteratively evaluating the Karush-Kuhn-Tucker (KKT) conditions to determine the contact state of any non-mortar element. The procedure stops when the partition into stick, slip and open regions remains unchanged after two consecutive active-set state update.
At each active-set iteration $l$, the inequality constraint \eqref{eq:constraint_inequality} can be replaced by the equality
\begin{equation}
    \mathcal{R}_t = \langle \boldsymbol \mu, \boldsymbol g_h \rangle_{\Gamma_f^{\text{stick}}} + \langle \mu_N, g_{N,h} \rangle_{\Gamma_f^{\text{slip}}} + \frac{1}{c}\langle \boldsymbol \mu_T, \boldsymbol t_{T,h} - \boldsymbol t_T^* \rangle_{\Gamma_f^{\text{slip}}} + \frac{1}{c}\langle \boldsymbol \mu, \boldsymbol t_h \rangle_{\Gamma_f^{\text{open}}} = 0,
    \label{eq:ActiveSet_residual}
\end{equation}
where $c$ is a unitary coefficient that restores the dimensional consistency.
To write the problem in matrix form, let us introduce into Eqs. \eqref{eq:displacement_equation} and \eqref{eq:ActiveSet_residual} the finite element approximations:
\begin{equation}
    \boldsymbol{u}^{(i)}_h = \sum_{j=1}^{|\mathcal{I}_u^{(i)}|} \boldsymbol{\eta}^{(i)}_j u_j, \qquad \boldsymbol{t}_h = \sum_{q=1}^{|\mathcal{I}_t|} \boldsymbol{\mu}_q t_q,
\end{equation}
then collect the displacement and traction discrete variables into vectors $\disp \in \mathbb{R}^{\card{\mathcal{I}_u}}$ and $\trac \in \mathbb{R}^{\card{\mathcal{I}_t}}$. 
By differentiating $\mathcal{R}_u$ and $\mathcal{R}_t$ with respect to $\disp$ and $\trac$, respectively, we obtain the Jacobian matrix:
\begin{equation}
    J = \begin{bmatrix}
        \frac{\partial \mathcal{R}_u}{\partial \disp} & \frac{\partial \mathcal{R}_u}{\partial \trac} \\
        \frac{\partial \mathcal{R}_t}{\partial \disp} & \frac{\partial \mathcal{R}_t}{\partial \trac}
    \end{bmatrix} =
    \begin{bmatrix}
        A & B_1 \\ B_2 & C 
    \end{bmatrix}
\end{equation}
with entries:
\begin{align}
    [A]_{jq} &= \sum_i(\nabla^s\boldsymbol{\eta}^{(i)}_j,\mathbb{D} : \nabla^s\boldsymbol {\eta}^{(i)}_k)_{\Omega^{(i)}}, \quad j,q=1,\ldots,|\mathcal{I}_u|, \nonumber\\
    [B_1]_{jq} &= \langle \llbracket \boldsymbol \eta_j \rrbracket, \boldsymbol \mu_q\rangle_{\Gamma_f^{(1)}}, \quad j=1,\ldots,|\mathcal{I}_u|, \; q=1,\ldots,|\mathcal{I}_t| \nonumber\\
    [B_2]_{jq} &= \langle \boldsymbol \mu_j, \jump{\boldsymbol\eta_q}\rangle_{\Gamma_f^{\text{stick}}} + \langle \mu_{j,N}, \llbracket \eta_{q,N} \rrbracket \rangle_{\Gamma_f^{\text{slip}}} - \frac{1}{c} \bigg\langle \boldsymbol{\mu}_{j,T}, \bigg( \frac{\partial \boldsymbol{t}_T^*}{\partial \Delta_n\boldsymbol{g}_T}\bigg)\llbracket \boldsymbol{\eta}_{q,T} \rrbracket \bigg\rangle_{\Gamma_f^{\text{slip}}}, \nonumber \\
    & \qquad j=1,\ldots,|\mathcal{I}_t|, \; q=1,\ldots,|\mathcal{I}_u|, \nonumber \\
    [C]_{jq} &= \frac{1}{c}\langle \boldsymbol\mu_{j,T},\boldsymbol\mu_{q,T} \rangle_{\Gamma_f^\text{slip}}  - \frac{1}{c} \bigg\langle \boldsymbol{\mu}_{j,T}, \bigg( \frac{\partial \boldsymbol{t}_T^*}{\partial t_N}\bigg) \mu_{q,N}\bigg\rangle_{\Gamma_f^{\text{slip}}}  + \frac{1}{c} \langle \boldsymbol \mu_j,\boldsymbol \mu_q \rangle_{\Gamma_f^{\text{open}}}, \nonumber \\
    & \qquad j,q=1,\ldots,|\mathcal{I}_t|. \nonumber
\end{align}
The non-linearity stems from the limiting tangential traction $\boldsymbol t_T^*$, whose derivatives with respect to the primary variables read:
\begin{align}
         \frac{\partial \boldsymbol{t}_T^*}{\partial \Delta_n\gt} &= \tau_{lim}(t_N)\frac{\|\Delta_n\gt\|^2\mathbf{I}-\Delta_n\gt \otimes \Delta_n\gt}{\| \Delta_n\gt\|^3}, \nonumber \\
     \frac{\partial \boldsymbol{t}_T^*}{\partial t_N} &= -\tan(\phi)\frac{\Delta_n\gt}{\| \Delta_n\gt \|}, \nonumber
\end{align}
By denoting with $\mathbf r_u$ and $\mathbf r_t$ the discrete residual contributions obtained from $\mathcal{R}_u$ and $\mathcal{R}_t$,
at each Newton iteration $k$ the solution is determined by 
\begin{align}
    &\text{solving} \ \ \  \begin{bmatrix}
        A & B_1 \\ B_2 & C 
    \end{bmatrix}^{l,(k)} \begin{bmatrix}
        \delta \disp \\ \delta \trac
    \end{bmatrix} = - \begin{bmatrix}
        \mathbf{r}_u \\ \mathbf r_t
    \end{bmatrix}^{l,(k)},   \label{eq:activeSet_system} \\ 
    &\text{setting} \ \ \ \begin{bmatrix}
        \disp \\ \trac
    \end{bmatrix}^{l,(k+1)} = \begin{bmatrix}
        \disp \\ \trac
    \end{bmatrix}^{l,(k)} + \begin{bmatrix}
        \delta \disp \\ \delta \trac
    \end{bmatrix}. \nonumber
\end{align}

Notice that, in the computation of the entries of $B_1$ and $B_2$, we have to consider the jump of $\boldsymbol{\eta}$, which involves the projection of $\boldsymbol{\eta}^{(2)}$ onto $\Gamma_f^{(1)}$. This generates a set of well-known issues in the practical implementation of the mortar method, which requires the integration of product of functions defined on different, non-conforming domains, which are often discontinuous on $\Gamma_f^{(1)}$.
For example, the entries of the contribution $\langle \boldsymbol{\mu}_j,\jump{\boldsymbol{\eta}_q} \rangle_{\Gamma_f^{\text{stick}}}$ into $[B_2]_{jq}$ read:
\begin{equation}
    \langle \boldsymbol{\mu}_j,\jump{\boldsymbol{\eta}_q} \rangle_{\Gamma_f^{\text{stick}}} = \langle \boldsymbol \mu_j, \boldsymbol \eta_q^{(1)}\rangle_{\Gamma_f^{\text{stick}}} - \langle \boldsymbol \mu_j, \Pi \boldsymbol \eta_q^{(2)}\rangle_{\Gamma_f^{\text{stick}}}.
    \label{eq:mortat_mat}
\end{equation}
The first inner product on the right-hand side of \eqref{eq:mortat_mat} gives rise to a standard mass matrix 
\cite{popp2010dual,yang2005two}. By distinction, the last contribution at the right-hand side of \eqref{eq:mortat_mat} involves the operator $\Pi$, which maps basis functions from $\Gamma_f^{(2)}$ to $\Gamma_f^{(1)}$. $\Pi$ is typically an orthogonal projector along an averaged normal computed on the non-mortar side.  
Integrating the product of basis functions having support on unrelated grids requires the adoption of non-standard quadrature algorithms and stands as the major algorithmic challenge of the mortar method \cite{maday2002influence,farah2015segment,moretto2024novel}. 

\section{Stability analysis}
The stability analysis of problems \eqref{eq:displacement_equation} and \eqref{eq:ActiveSet_residual} becomes critical when some fracture elements are in the stick regime. For simplicity, we assume that $\Gamma_f^{(1)} \equiv \Gamma_f^{\text{stick}}$, i.e., the entire contact surface is in the stick condition. In this case, system \eqref{eq:activeSet_system} becomes
\begin{align}
\begin{bmatrix}
        A & B \\ B^T & 0 
    \end{bmatrix} \begin{bmatrix}
        \delta \disp \\ \delta \trac
    \end{bmatrix} = - \begin{bmatrix}
        \mathbf{r}_u \\ \mathbf r_t
    \end{bmatrix},    
    \label{eq:saddlePointSystem}
\end{align}
where, for the sake of clarity, we set $B\equiv B_1$ and suppressed the unused superscripts $l$ and $k$, referring to the active-set and Newton iteration counts, respectively, since all arguments are valid for any active-set and nonlinear iteration.
The system \eqref{eq:saddlePointSystem} recovers a saddle point problem that must satisfy the inf-sup condition
\begin{align}
    \inf_{\boldsymbol{\mu}_h \in \boldsymbol{\mathcal{M}}_h} \sup_{\boldsymbol{v}_h \in \boldsymbol{\mathcal{V}}_h} \frac{\langle \boldsymbol{\mu}_h, \boldsymbol{g}_h \rangle_{\Gamma_f^{(1)}}}{\|\boldsymbol \mu_h \|_{\boldsymbol{\mathcal M}} \|\boldsymbol v_h \|_{\boldsymbol{\mathcal V}}} \geq \beta > 0.
    \label{eq:inf_sup_condition}
\end{align}
for its solvability.

In the proposed formulation, we use the finite element spaces:
\begin{subequations} \label{eq:discrete_spaces}
\begin{align} \label{eq:femSpace_primaryVar}
    \boldsymbol{\mathcal{V}}_h = \{ \boldsymbol v_h \in \boldsymbol{\mathcal{V}}: \boldsymbol{v}_h \vert_{K} \in [\mathbb{Q}_1(K)]^3, \forall K  \in \mathcal{T} \}, \\
    \label{eq:femSpace_multipliers}
    \boldsymbol{\mathcal M}_h = \{ \boldsymbol \mu_h \in \boldsymbol{\mathcal{M}}: \boldsymbol{\mu}_h \vert_{F} \in [\mathbb{P}_0(F)]^3, \forall F  \in \mathcal{F} \}, 
\end{align}    
\end{subequations}
where $[\mathbb{Q}_1]^3$ is the space of trilinear nodal basis functions for hexahedral elements, while $[\mathbb{P}_0]^3$ is the space of piecewise constant basis functions.
It is well known that this pair of spaces does not satisfy the inf-sup condition \eqref{eq:inf_sup_condition}. However, the use of $\mathbb{P}_0$ offers several beneficial properties for the implementation of the mortar method. For instance, Dirichlet boundary conditions can be applied on the fault boundary without local modification of the multipliers basis functions \cite{puso2003mesh}. Similarly, solvability issues due to cross-points of intersecting interfaces are inherently addressed.
These features are particularly important in geological applications, where intricate geometries with numerous intersections and complex boundary conditions are frequently encountered.
Therefore, we have to introduce a stabilization strategy to satisfy \eqref{eq:inf_sup_condition} while keeping the spaces defined in \eqref{eq:discrete_spaces}.

Inspired by the pressure-jump stabilization \cite{elman2014finite} used in the Stokes community, we develop a traction-jump stabilization by considering the stabilized constraint equation
\begin{align}
    \mathcal{R}_t^* = \mathcal{R}_t + s(\boldsymbol{\mu},\boldsymbol t_h),  
\end{align}
where $s$ is a bilinear form introducing the stabilizing contribution:
\begin{align}
    s(\boldsymbol \mu,\boldsymbol t_h)  = \sum_{E_{KL} \in \mathcal{E}} \jump{\boldsymbol{\mu}}_{E_{KL}} \mathbf{S}_{E_{KL}} \jump{\boldsymbol t_h}_{E_{KL}}. 
    \label{eq:stabilizationTerm}
\end{align}
In equation \eqref{eq:stabilizationTerm}, $\jump{\bullet}_{E_{KL}}$ is the jump operator across the edge $E_{KL}$, defined as $(\bullet) \vert_{F_K} - (\bullet) \vert_{F_L}$, and $\mathbf{S}_{E_{KL}}$ is a second order tensor that provides the appropriate scaling to the local stabilization contribution. 
A detailed analysis of the inf-sup stability of the proposed formulation can be found in \cite{burman2014projection}, where the jump stabilization has been recast into the framework of projection-stabilization techniques. 

\subsection{Stabilization derivation}
\begin{figure}
    \centering
    \includegraphics[width=0.6\linewidth]{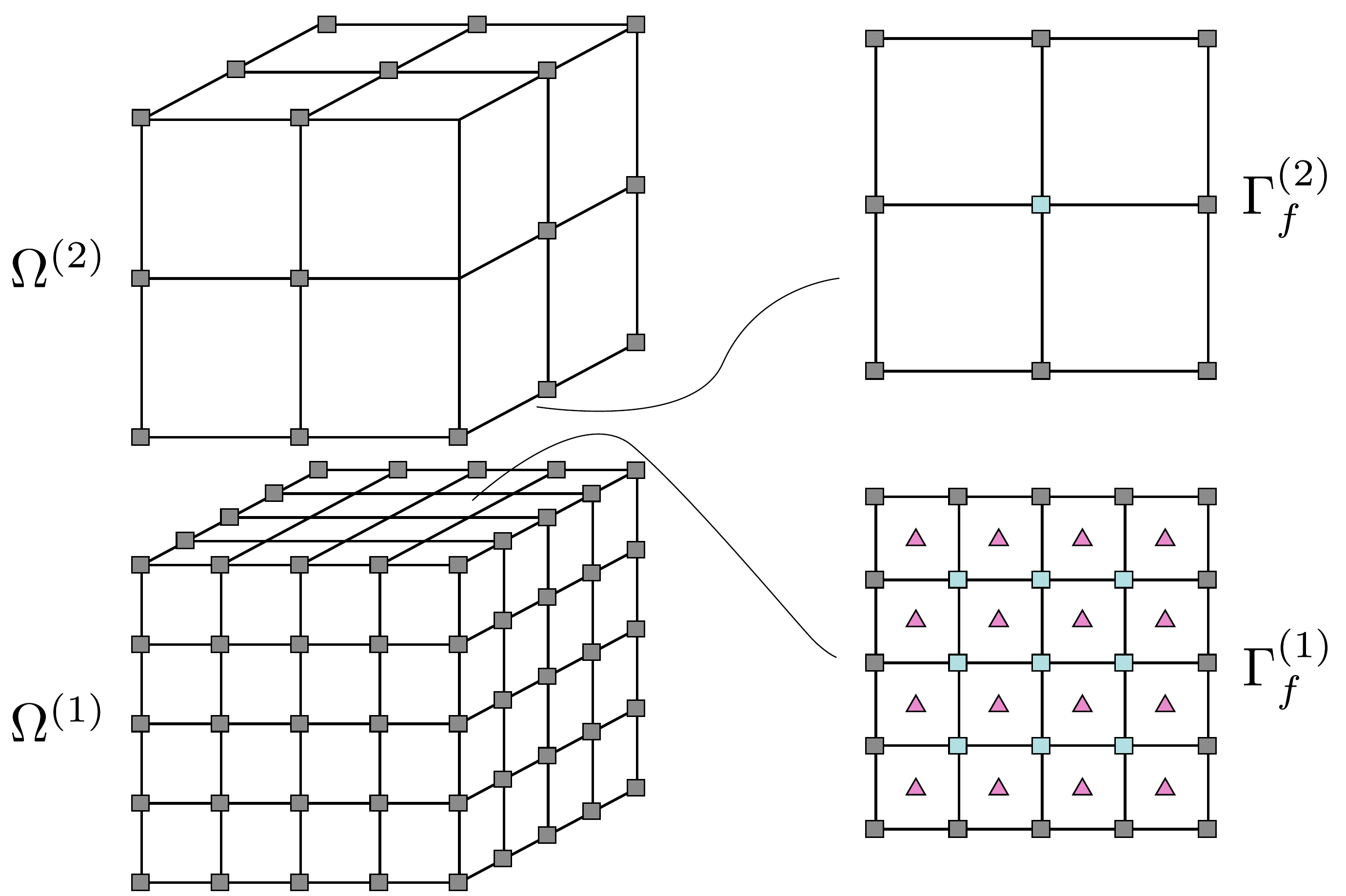}
    \caption{Base Macroelement for the instability analysis and the inf-sup constant evaluation ($r_h = 1/2)$, with location of fixed displacement \filledsquare[gray]{1pt}, free displacements \filledsquare[CornflowerBlue]{1pt} and tractions \filledtriangle[magenta]{8pt}.}
    \label{fig:basePatch}
\end{figure}

We now derive an approach to evaluate the scaling tensor $\mathbf{S}_{E_{KL}}$, building upon the macroelement approach proposed in \cite{franceschini2020algebraically} for the conforming case and extending it to the non-conforming one. 
We first identify the origin of the instability by considering an unstable non-conforming patch of elements (Figure \ref{fig:basePatch}), where $r_h = h^{(1)}/h^{(2)} = 1/2$ and all the displacement degrees of freedom are fixed at the boundary. By assembling the Jacobian matrix in \eqref{eq:saddlePointSystem}, one can observe that its Schur complement $S = B^T A^{-1}B$ has a nontrivial kernel, i.e., $\ker(S) \neq \emptyset$. The eigenmodes of the kernel of the Schur complement correspond to spurious traction profiles and reveal the instability of the formulation.

For the stabilization term \eqref{eq:stabilizationTerm} to be effective, we have to ensure that the resulting stabilized problem,
\begin{align}
\begin{bmatrix}
        A & B \\ B^T & -H 
    \end{bmatrix} \begin{bmatrix}
        \delta \disp \\ \delta \trac
    \end{bmatrix} = - \begin{bmatrix}
        \mathbf{r}_u \\ \mathbf r_t
    \end{bmatrix},    
    \label{eq:stabilizedSaddlePointSystem}    
\end{align}
is such that 
\begin{align}
    \ker(S+H) = \emptyset.
    \label{eq:stabilizedKernel}
\end{align} 
As proved in \cite{moretto2025stabilized}, the jump stabilization term satisfies this relationship by construction for every positive definite $\mathbf{S}_{E_{KL}}$ tensor. However, a proper tuning of the entries in $\mathbf{S}_{E_{KL}}$ is needed to ensure that the stabilization contribution is not too small to effectively remove the spurious modes, and not too large to avoid an excessive loss of accuracy in the primary variable due to overly smooth tractions. In fact, the stabilizing contribution modifies the residual equation $\mathcal{R}_t$ enforcing the constraint. Therefore, an unnecessarily large stabilization term will cause an excessive deviation from the original constraint.

In view of \eqref{eq:stabilizedKernel}, it is natural to scale $H$ so that its entries have the same order as those of the Schur complement $S$. Therefore, we propose here an implementation strategy to compute $\mathbf{S}_{E_{KL}}$ as a local Schur complement approximation. 
We consider the macroelement shown in Figure \ref{fig:basePatch}, consisting of mortar cells connected to an internal interface node (the blue one in the center of $\Gamma_f^{(2)}$) and the set of non-mortar cells sharing some supports with the mortar cells. We use the superscript $\check{\bullet}$ to denote any set or quantity restricted to the macroelement. 
Once the macroelement has been defined, we can compute an approximation of the local Schur complement as
\begin{align}
    \tilde{S} = \check{B}\check{D}_A^{-1}\check{B}{^T},
\end{align}
where $\check{D}_A = \operatorname{diag}(\check{A})$. We keep only the diagonal entries of the local gathering from the stiffness matrix, 
as this provides the correct order of magnitude and allows a straightforward inversion.
Then, for each internal non-mortar edge in the macroelement, we assemble the jump stabilization contribution \eqref{eq:stabilizationTerm} with $\mathbf{S}_{E_{KL}}$ defined as
\begin{align}
    \mathbf{S}_{E_{KL}} = \frac{1}{2}(\tilde{S}_{F_K} +  \tilde{S}_{F_L}),
\label{eq:localStabilizationScaling}
\end{align}
where $\tilde{S}_{F_i}$ is the submatrix containing the entries related to the traction degrees of freedom of the face $F_i$.
We repeat this computation for all the internal mortar nodes and their related macroelement. 

\begin{figure}
    \centering
    \includegraphics[width=0.7\linewidth]{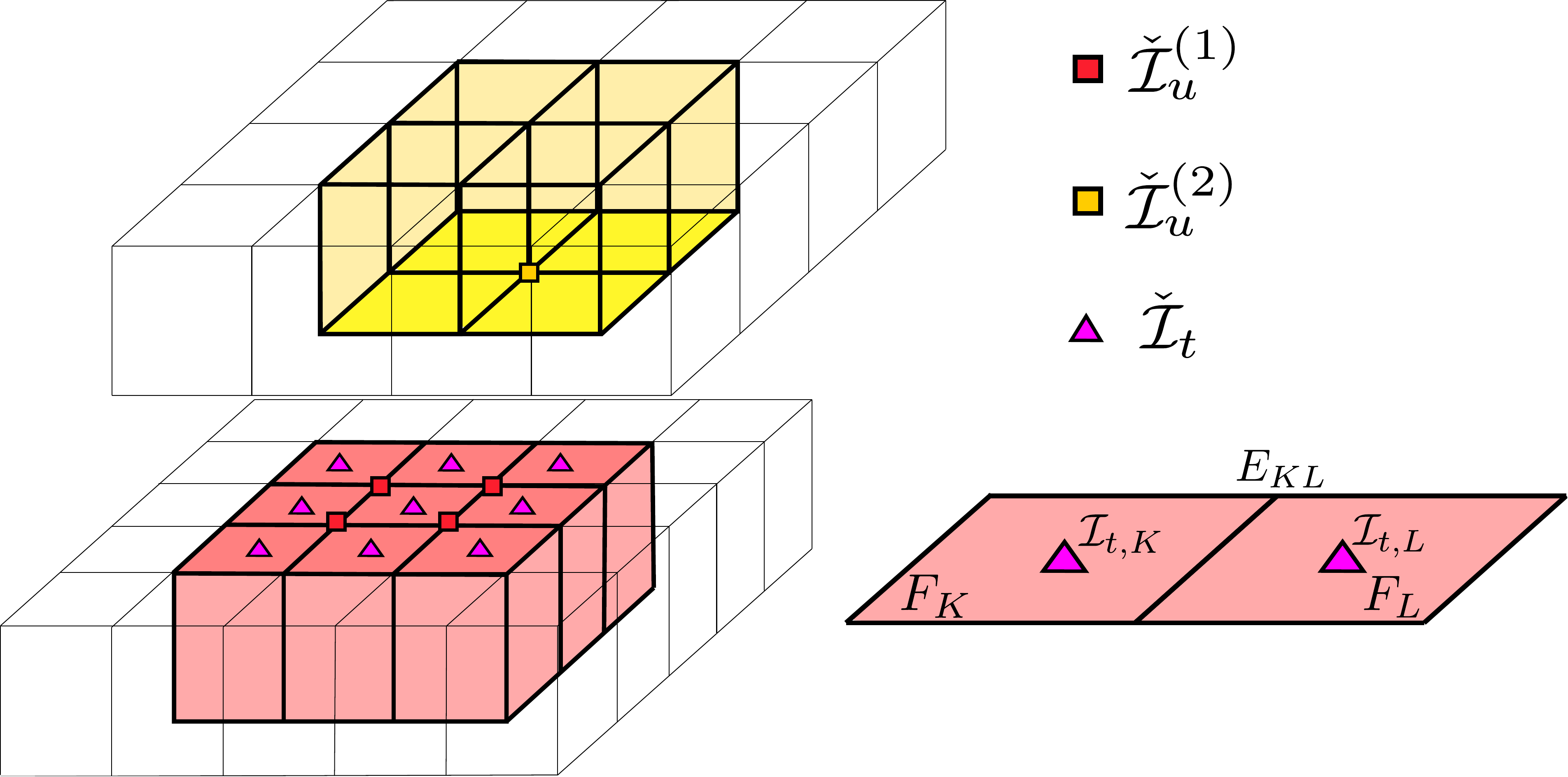}
    \caption{Local macroelement for the computation of $\mathbf{S}_{E_{KL}}$.}
    \label{fig:localMacroelement}
\end{figure}

Algorithm \ref{alg:stabMatrix} summarizes the procedure using a MATLAB-style notation with reference to the macro-element and the related symbols in Figure \ref{fig:localMacroelement}.
The entries of $\check{B}$ are computed using the mortar matrices $D$ and $M$ corresponding to the discretization of the inner products introduced in \eqref{eq:mortat_mat}:
\begin{equation}
[D]_{jq} = \langle \boldsymbol \mu_j, \boldsymbol \eta_q^{(1)}\rangle_{\Gamma_f^{(1)}}, \quad \text{and} \quad
[M]_{jq} = \langle \boldsymbol \mu_j, \Pi \boldsymbol \eta_q^{(2)}\rangle_{\Gamma_f^{(1)}}. \nonumber
\end{equation}

\begin{algorithm}
\caption{Computation of traction-jump stabilization matrix}
    \begin{algorithmic}[1]
    \For{\textbf{all} mortar nodes}
    \State Build macroelement and find $\check{\mathcal{E}}$, $\check{\mathcal{I}}_{u}^{(1)},\check{\mathcal{I}}_{u}^{(2)},\check{\mathcal{I}}_{t}$ 
    \State Gather local $\check A=\texttt{blkdiag}(A(\check{\mathcal{I}}_{u}^{(1)},\check{\mathcal{I}}_{u}^{(1)}),A(\check{\mathcal{I}}_{u}^{(2)},\check{\mathcal{I}}_{u}^{(2)}))$
    \State Compute $\check D_{A} = \texttt{diag}(\check A)$
    \State Gather local $\check B^T = [D(\check{\mathcal{I}}_{t},\check{\mathcal{I}}_{u}^{(1)}),-M(\check{\mathcal{I}}_{t},\check{\mathcal{I}}_{u}^{(2)})]$ 
    \State Compute $\tilde{S} = \check B \check D^{-1}_{A} \check B^T$
    \For{$E_{KL}  \text{ in } \check{\mathcal{E}}$}
    \State Compute $\mathbf{S}_{E_{KL}}$ with \eqref{eq:localStabilizationScaling} 
    \State Compute local $\check H = [\mathbf{S}_{E_{KL}}, -\mathbf{S}_{E_{KL}};-\mathbf{S}_{E_{KL}}, \mathbf{S}_{E_{KL}}]$
    \State Find $\mathcal{I}_{t,K},\mathcal{I}_{t,L} \in \check{\mathcal{I}}_{t}$
    \State Assemble $H([\mathcal{I}_{t,K},\mathcal{I}_{t,L}],[\mathcal{I}_{t,K},\mathcal{I}_{t,L}]) = H([\mathcal{I}_{t,K},\mathcal{I}_{t,L}],[\mathcal{I}_{t,K},\mathcal{I}_{t,L}])+\check H$ 
    \EndFor
    \EndFor
    \end{algorithmic}
    \label{alg:stabMatrix}
\end{algorithm}

\subsection{Stabilization verification}

We verify the effectiveness of the proposed stabilization by numerically checking the inf-sup condition \eqref{eq:inf_sup_condition}. We follow the procedure described in \cite{el2001stability,elman2014finite}, which is briefly outlined here for convenience.
First, we use the bound 
\begin{align}
    \| \boldsymbol \mu_h \|_{\boldsymbol{\mathcal{M}}} = \|\boldsymbol \mu_h\|_{-\frac{1}{2};\Gamma_f^{(1)}} \geq ch^{\frac{1}{2}}  \|\boldsymbol \mu_h\|_{0;\Gamma_f^{(1)}}, \ \ \ \ \forall \, \boldsymbol{\mu}_h \in \boldsymbol{\mathcal{M}}_h
    \label{eq:normEstimate}
\end{align}
where $c$ is a positive constant. We can use \eqref{eq:normEstimate} to reformulate \eqref{eq:inf_sup_condition} into the weaker condition
\begin{align}
    \inf_{\boldsymbol \mu_h \in \boldsymbol{\mathcal{M}}_h} \sup_{\boldsymbol v_h \in \boldsymbol{\mathcal{V}}_h} \frac{\langle \boldsymbol \mu_h, \boldsymbol v_h\rangle_{\Gamma_f^{(1)}}}{h^{\frac{1}{2}}\|\boldsymbol \mu_h \|_{0;\Gamma_f^{(1)}} \|\boldsymbol v_h\|_{\boldsymbol{\mathcal{V}}}} = \beta^*  \geq \beta > 0,
    \label{eq:inf-sup_weaker}
\end{align}
in order to avoid the evaluation of the fractional dual norm of the multipliers. 
However, a test based on \eqref{eq:inf-sup_weaker} provides only a necessary condition to satisfy the inf-sup condition. 
Hence, we rewrite \eqref{eq:inf-sup_weaker} in a discrete form as
\begin{align}
    \min_{\trac} \max_{\disp} \frac{\trac^TB^T\disp}{h^{\frac{1}{2}}\sqrt{\trac^T Q \boldsymbol \trac} \sqrt{\disp^T A \disp }} = \beta_h^* > 0,
    \label{eq:inf-sup_matrix}
\end{align}
after using Poincar\'e inequality for the displacement norm and having introduced the matrix $Q$ such that $\mathbf{t}^T Q\mathbf{t} = \langle \boldsymbol t_h, \boldsymbol t_h\rangle_{\Gamma_f^{(1)}}$.
Following the derivation of \cite{elman1996iterative,elman2014finite}, we obtain from \eqref{eq:inf-sup_matrix} the following characterization of the inf-sup constant:
\begin{align}
    (\beta^*)^2 = \min_{\trac \neq \boldsymbol{0}} \frac{1}{h} \frac{\trac^T (B A^{-1} B^T) \trac}{\trac^T Q \trac},
\end{align}
that is, $\beta^*$ is the square root of the smallest eigenvalue of $(hQ)^{-1}(B^TA^{-1}B)$, which corresponds to the scaled Schur complement of \eqref{eq:saddlePointSystem}.

\begin{figure}
    \centering
    \begin{subfigure}[b]{0.43\textwidth}
        \centering
        \includegraphics[width=\textwidth]{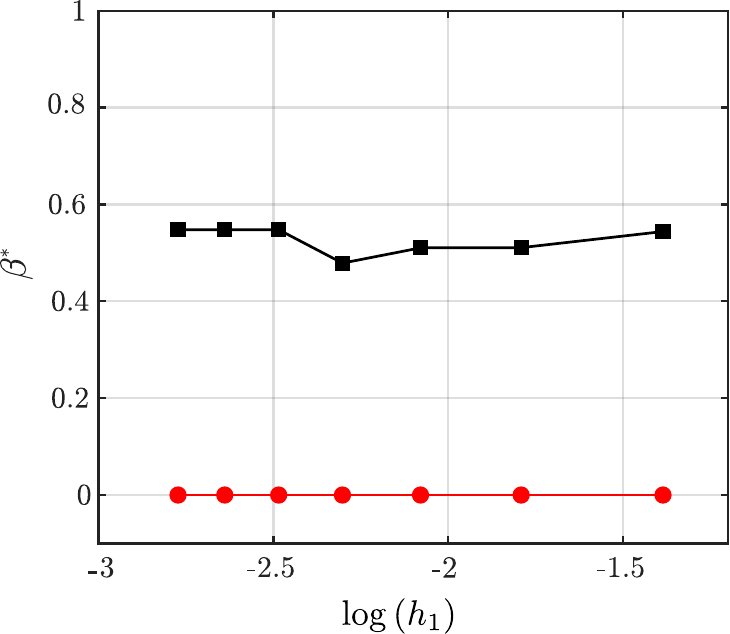}
        \caption{}
        \label{fig:infSupA}
    \end{subfigure}
    \hspace{0.1\textwidth}
    \begin{subfigure}[b]{0.4\textwidth}
        \centering
        \includegraphics[width=\textwidth]{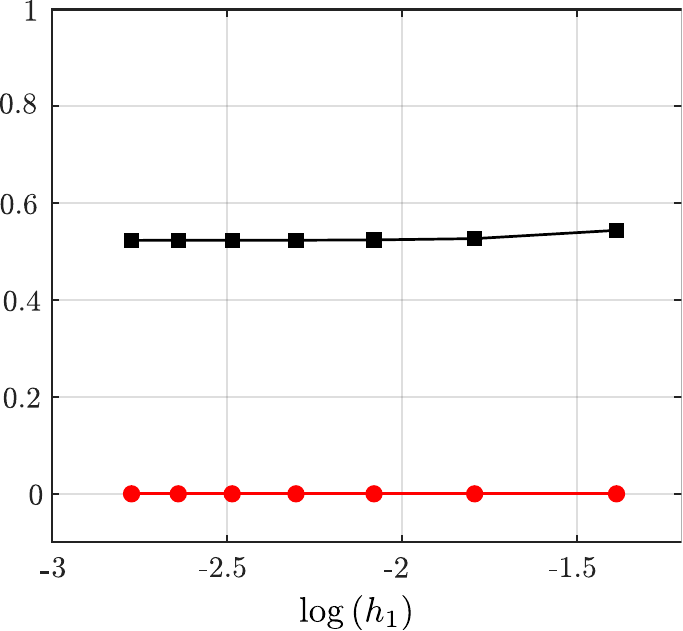}
        \caption{}
        \label{fig:infSupB}
    \end{subfigure}
    \caption{Results of the inf-sup test: refinement with constant $r_h$ (a), refinement on the non-mortar side only (b). Results with stabilized formulation in black line, unstabilized case in red line.}
    \label{fig:infSup}
\end{figure}

To perform the inf-sup test, we compute $\beta^*$ for subsequent refinements of the patch of Figure \ref{fig:basePatch}, while maintaining a fixed ratio $r_h$ between the mortar and non-mortar mesh sizes. The test is passed if the inf-sup constant does not decrease with the mesh size. 
The results of the test are shown in Figure~\ref{fig:infSupA} for both the stabilized and the unstabilized case. 
As expected, $\beta=0$ for the unstabilized case, since the Schur complement has a kernel. 
In the stabilized case, the inf-sup constant remains bounded away from zero for all considered mesh refinements.

We repeat the test, now fixing the mortar mesh size and refining only the non-mortar side, in order to check the sensitivity of $\beta^*$ to the refinement ratio $r_h$. The test is passed also in this case (Figure \ref{fig:infSupB}), being $\beta^*$ bounded away from zero and almost constant.

\section{Numerical results}
In this section, we present several numerical examples to validate the contact algorithm and demonstrate the effectiveness and the potential of the proposed stabilization.
\begin{enumerate}
    \item First, we test the method against analytical benchmarks to verify the correctness and accuracy of the proposed formulation (Test 1, 2, and 3).
    \item Second, we propose a test case that shows the effectiveness of the traction-jump stabilization when the non-mortar side is finer than the mortar one (Test 4). The results are compared against LMM formulations that are theoretically inf-sup stable, but can show oscillations in this case. 
    \item Finally, we present a realistic problem dealing with a large discontinuity in geological media (Test 5). This test case demonstrates the applicability of the formulation to a problem representative of real-world subsurface applications discretized by a corner-point grid.
\end{enumerate}

\subsection{Test 1: Patch test}
\begin{figure}
    \centering
    \begin{subfigure}[b]{0.3\textwidth}
        \centering
        \includegraphics[width=\textwidth]{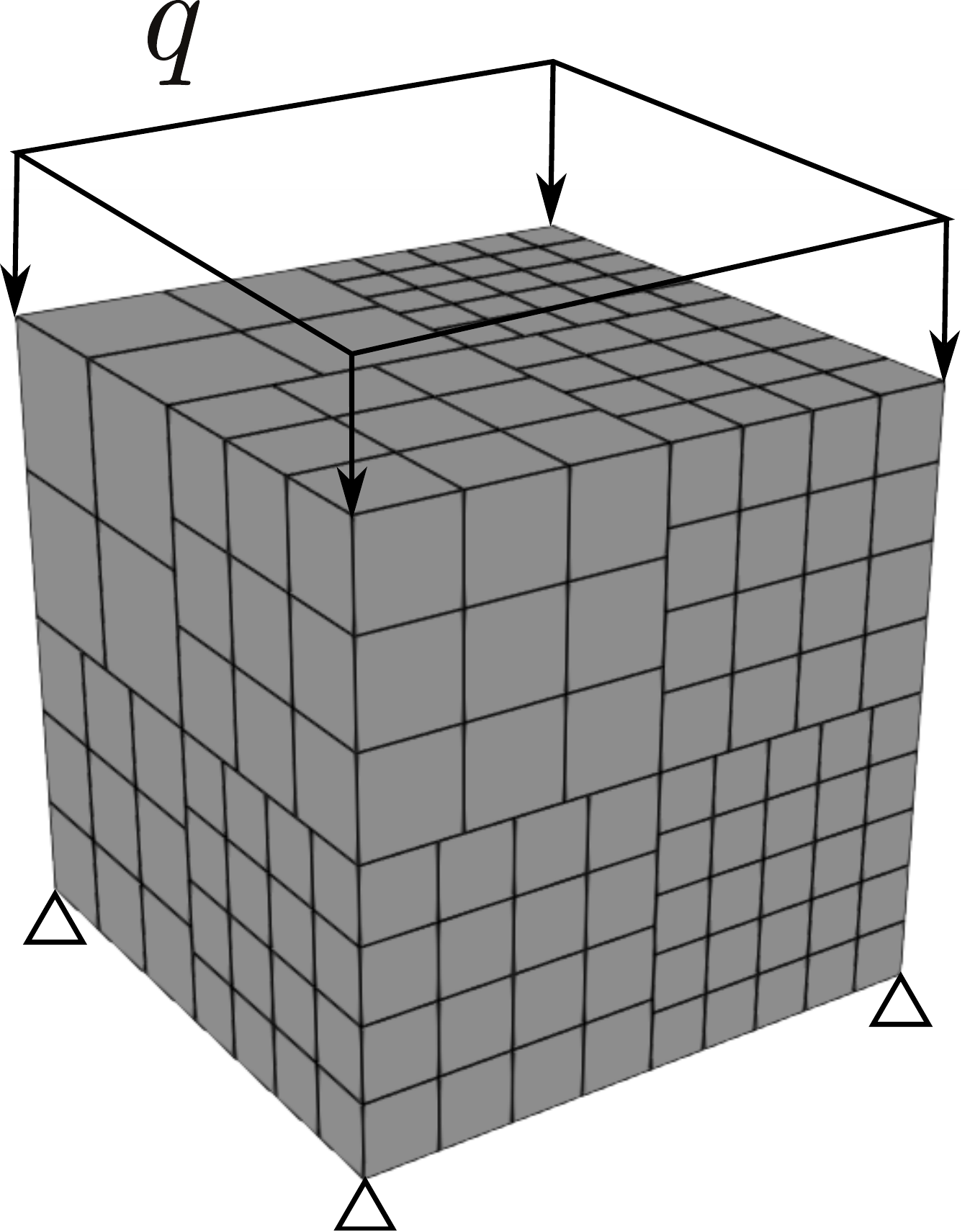}
        \caption{}
        \label{fig:PatchTest_setting}
    \end{subfigure}
    \hspace{0.15\textwidth}
    \begin{subfigure}[b]{0.45\textwidth}
        \centering
        \includegraphics[width=\textwidth]{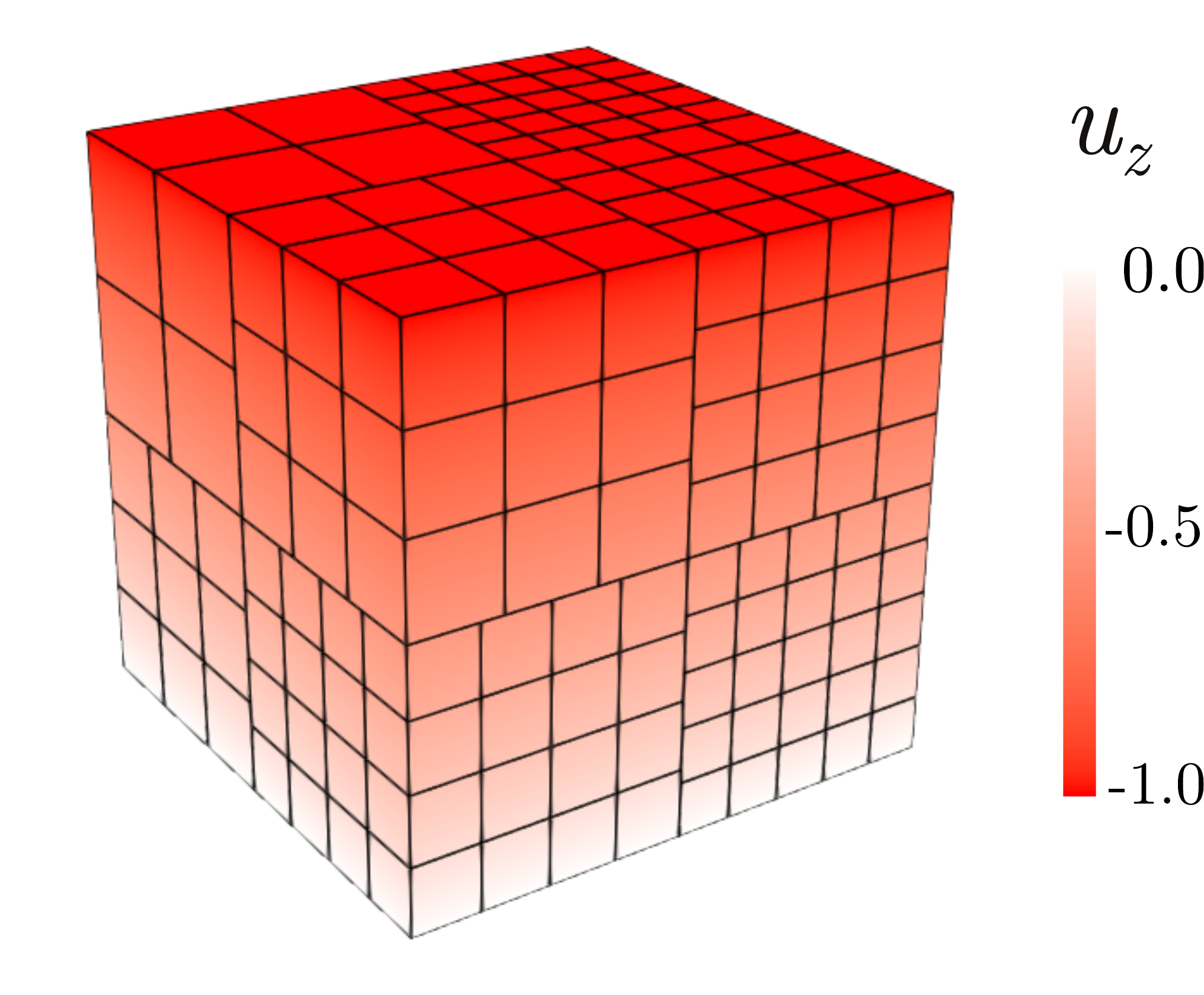}
        \caption{}
        \label{fig:PatchTest_uz}
    \end{subfigure}
    \caption{Test 1, patch test: model setting (a) and contour plot of the vertical displacement using the stabilized $\mathbb{P}_0$ formulation (b). }
    \label{fig:PatchTest}
\end{figure}

The first benchmark is a classical contact mechanics patch test, used to verify that the proposed stabilized formulation is capable of transferring a constant traction field across the interface. 
The domain is the unit cube, whose grid is made of eight blocks with different grid sizes, and intersecting at the cross point $(0.5,0.5,0.5)$ (see Figure \ref{fig:PatchTest_setting}). As discussed in the introduction, it is well known that the use of nodal Lagrange multipliers without an appropriate treatment of the basis functions at such cross points may lead to local oscillations, thereby violating the patch test. The use of stabilized $\mathbb{P}_0$ multipliers naturally addresses this issue.

\begin{figure}
    \centering
    \includegraphics[width=0.55\linewidth]{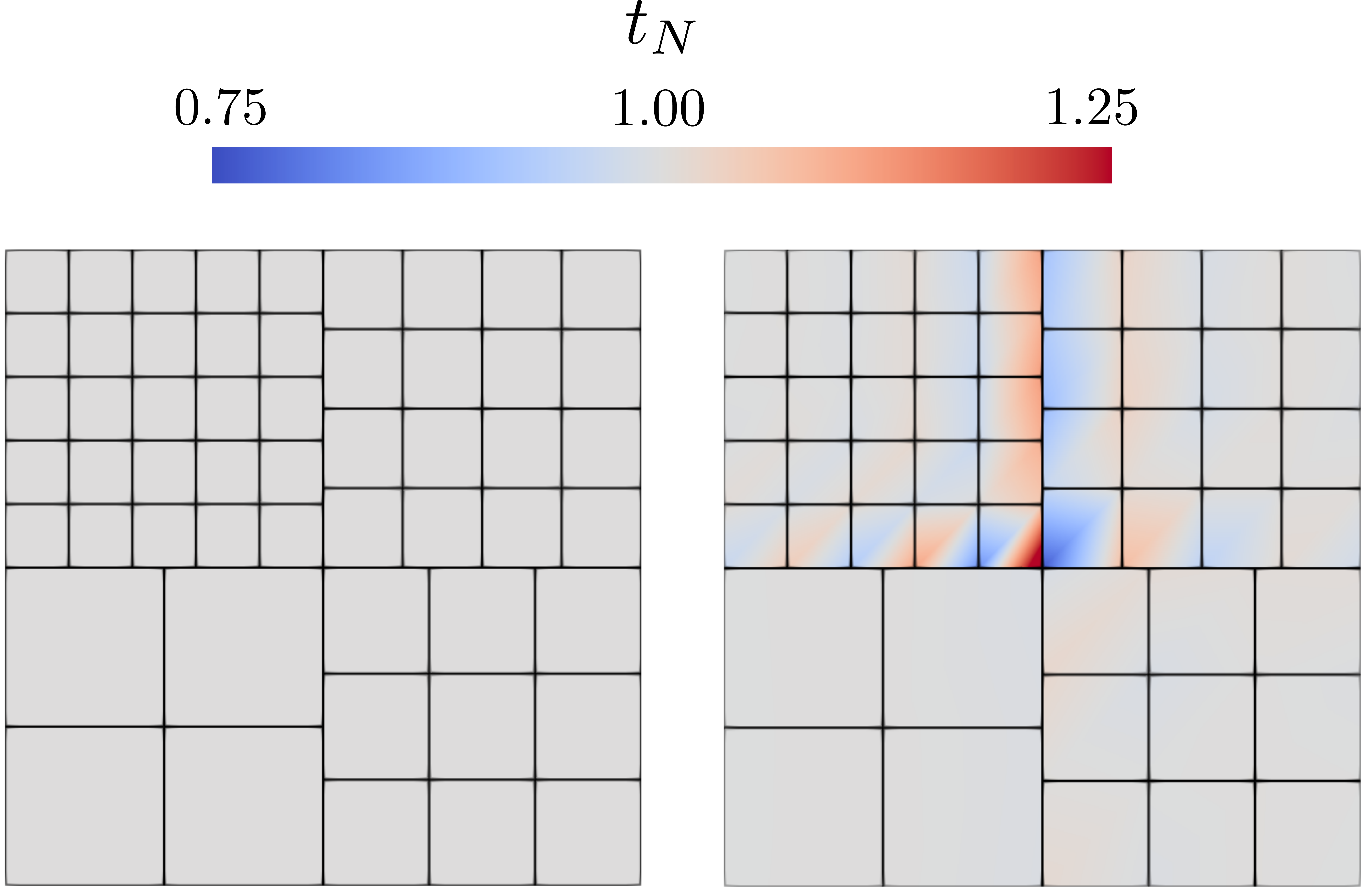}
    \caption{Test 1, patch test: normal stress at the interface for stabilized $\mathbb{P}_0$ multipliers (left) and standard nodal multipliers (right).}
    \label{fig:PatchTest_stress}
\end{figure}

The cube is vertically loaded by a constant traction $q = -1$ on the top face, and is fixed at the bottom face. We set the elastic properties as $E=1000$ and $\nu=0$ to prevent lateral deformation. The resulting vertical displacement field with $\mathbb{P}_0$ multipliers is reported in Figure \ref{fig:PatchTest_uz}.
For the sake of a comparison, we run the patch test also using a formulation with standard nodal multipliers, without any cross-point modification. 
Figure \ref{fig:PatchTest_stress} shows the normal traction computed at the interface. While the proposed formulation reproduces the exact constant $t_N=-1$, standard multipliers, although theoretically inf-sup stable, exhibit spurious oscillations along the non-conforming interfaces between different discretizations.

\subsection{Test 2: Constant sliding}
\begin{figure}
    \centering
    \begin{subfigure}[b]{0.24\textwidth}
        \centering
        \includegraphics[width=\textwidth]{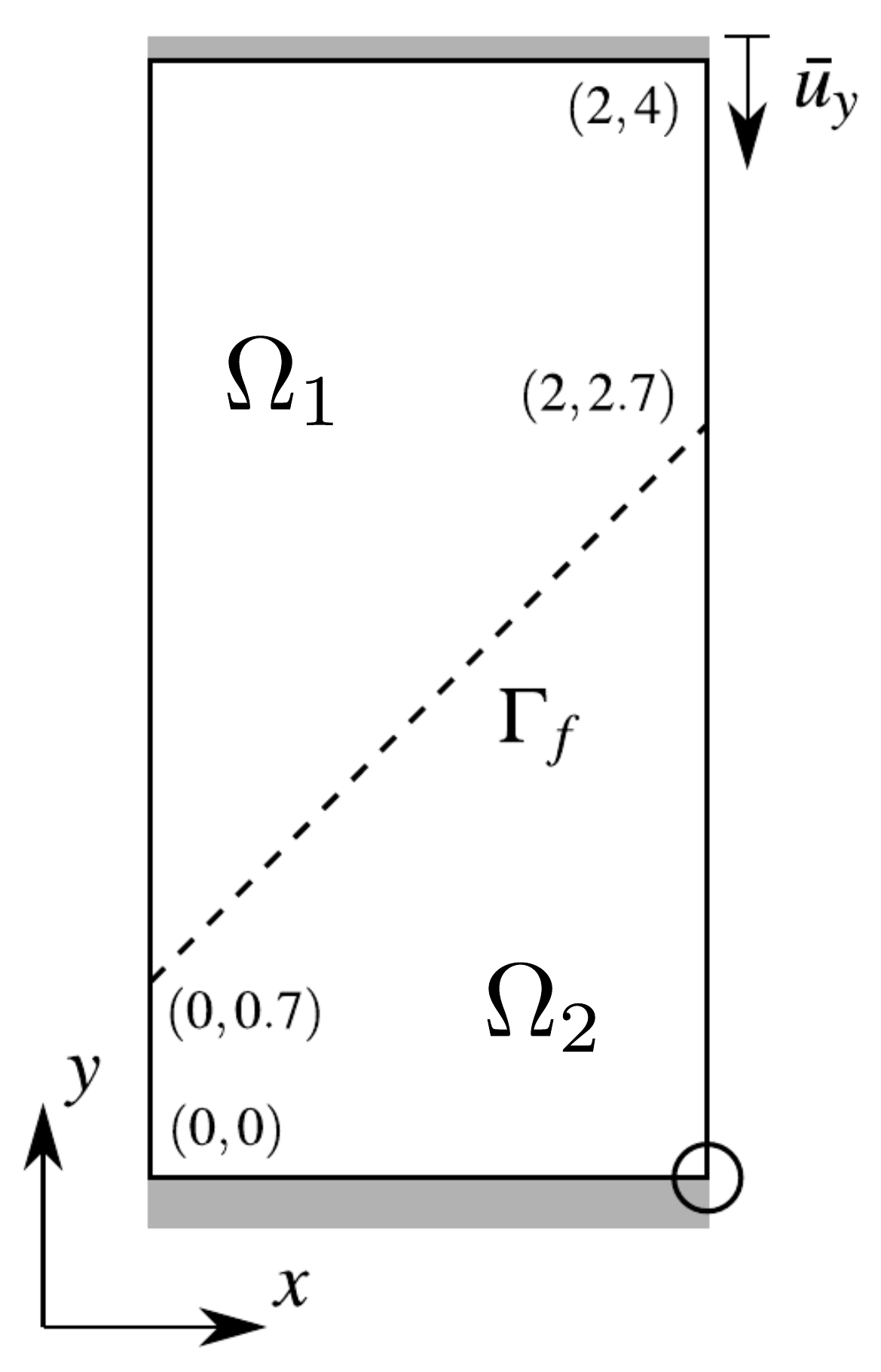}
        \caption{}
        \label{fig:constant_sliding_setting}
    \end{subfigure}
    \hspace{0.15\textwidth}
    \begin{subfigure}[b]{0.28\textwidth}
        \centering
        \includegraphics[width=\textwidth]{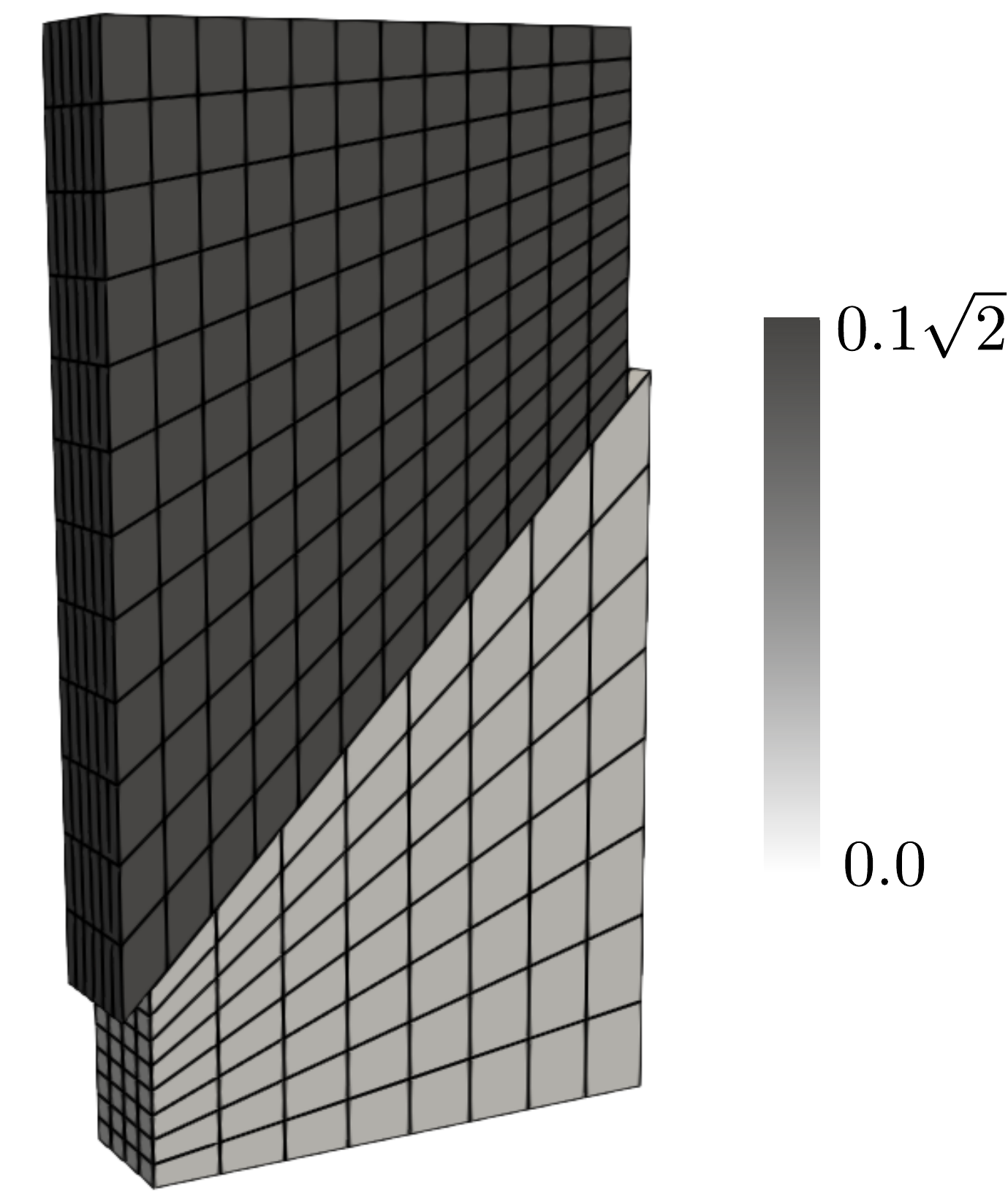}
        \caption{}
        \label{fig:constant_sliding_uz}
    \end{subfigure}
    \caption{Test 2, constant sliding: problem setting (a) and  contour of the displacement norm with deformed configuration (b).}
    \label{fig:constant_sliding}
\end{figure}

The second benchmark is a classical test case for sliding contact formulations. A two-dimensional representation of the domain is shown in Figure~\ref{fig:constant_sliding_setting}. 
A prismatic domain is divided by a tilted fault, with the non-mortar side at the bottom. 
A vertical displacement $\bar{u}_y$ is imposed on the top, while the vertical displacements are fixed at the bottom. The nodes along $(2,0,z)$ are fixed, with the circle in Figure~\ref{fig:constant_sliding_setting} showing their position.
Since the simulation is carried out in a fully 3D setting while the analytical solution is available only in 2D, a plane strain condition must be enforced by constraining the displacements along the $z$-axis. The material parameters are $E=5000$, $\nu=0.25$, $\phi=5.71^\circ$, $c=0$.
The analytical solution is a fully sliding interface with $\|\boldsymbol{g}_T \|_2= \tan{\phi}\,\sqrt{2}=0.1\sqrt{2}$. The proposed formulation exactly reproduces the analytical solution (Figure \ref{fig:constant_sliding_uz}), and requires only one active set update after the first elastic step.

\subsection{Test 3: Single fracture under compression}
\begin{figure}
    \centering
    \includegraphics[width=\linewidth]{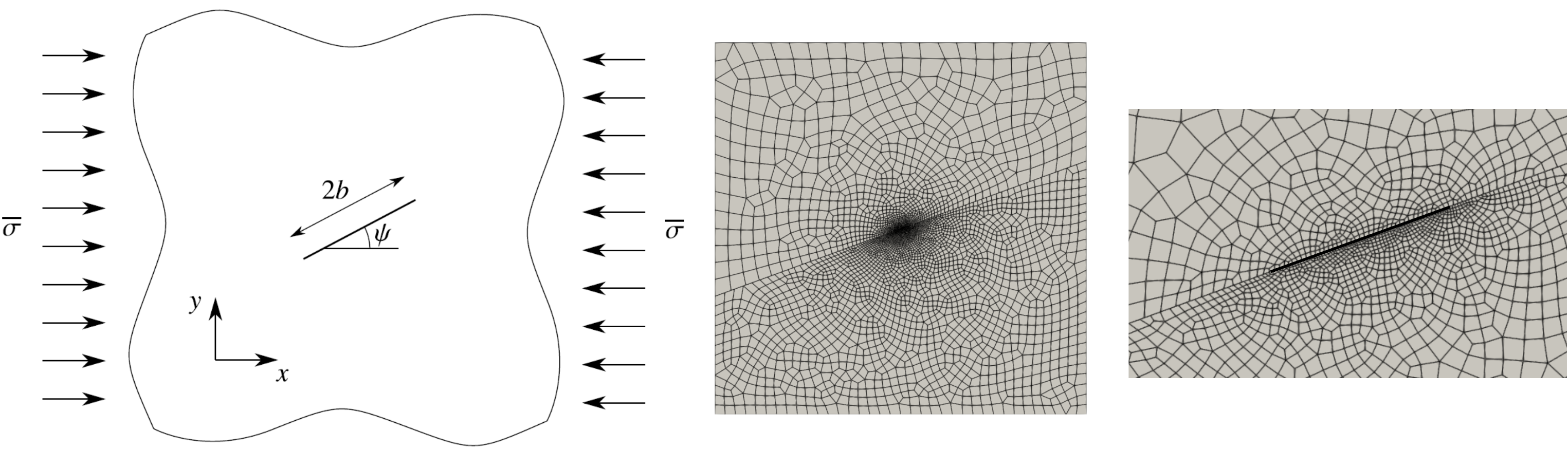}
    \caption{Test 3, single fracture under compression: problem setting (left), top view of the model grid (center) and zoom around the fracture (right).}
    \label{fig:SingleFracture_setting}
\end{figure}

The third benchmark considers a single crack embedded in a 2D infinite domain subjected to uniform uniaxial compression (Figure \ref{fig:SingleFracture_setting}). This reference problem is described in detail in \cite{phan2003symmetric}. Let $\bar{\sigma}=100$ be the magnitude of the applied compressive stress, $\psi=20^\circ$ the inclination angle of the fracture, \(2b=2\) its length, $\phi=30^\circ$ and zero cohesion. The analytical solution provides closed-form expressions for the normal traction and the tangential sliding along the fracture, given by
\begin{align}
t_N &= -\bar\sigma \sin^2\psi, \label{eq:normal_traction} \\[0.5em]
\|\boldsymbol{g}_T\| &=
\frac{4(1-\nu^2)}{E}
\left[ \bar\sigma \sin\psi(\cos\psi - \sin\psi\tan\phi)\right]
\sqrt{\, b^2 - (b-\xi)^2 \,}, \label{eq:sliding}
\end{align}
where $\xi$ denotes the curvilinear coordinate along the fracture, with $0 \le \xi \le 2b$. 
The problem is solved in a three-dimensional extruded $50\times50\times1$ domain, with appropriate plane strain conditions. We employ an unstructured grid with $r_h \approx 1/2$.
The elastic parameters are $E=15000$ and $\nu = 0.25$.

\begin{figure}
    \centering
    \includegraphics[width=0.9\linewidth]{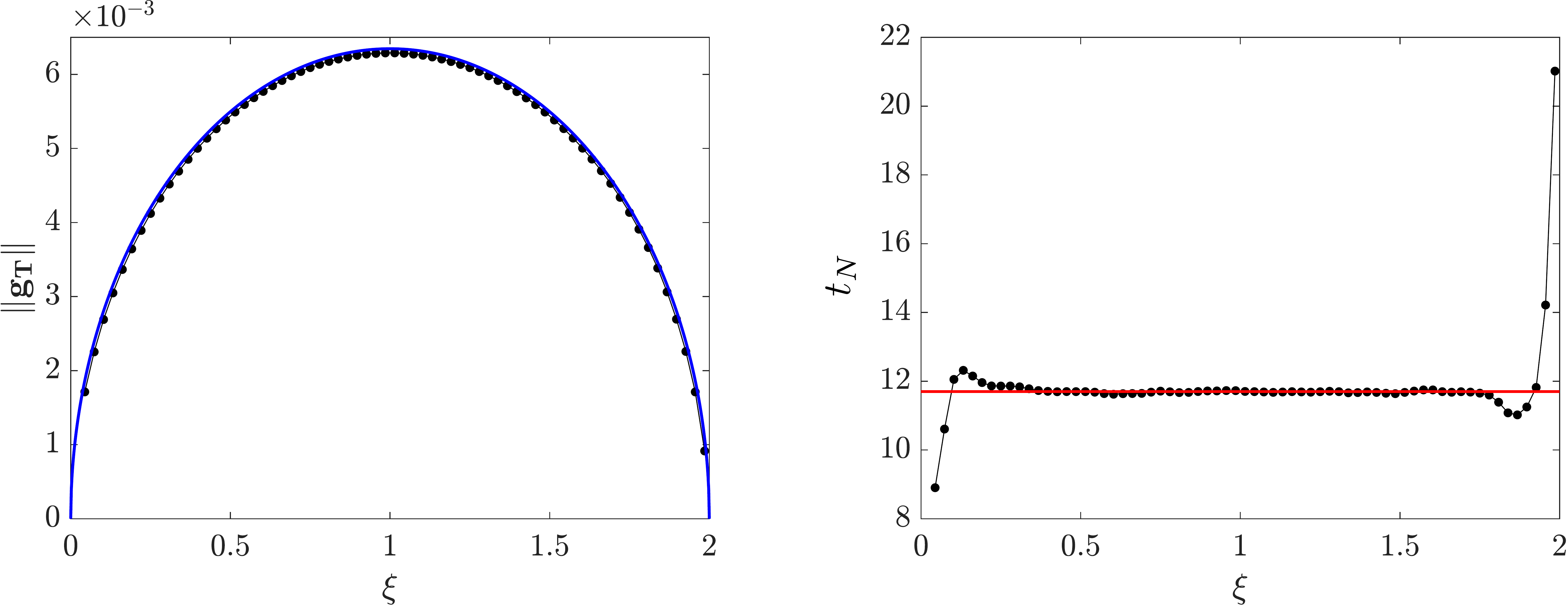}
    \caption{Test 3, single fracture under compression: $\|\boldsymbol{g}_T\|_2$ and $\sigma_N$ on the fracture. Analytical solution with continuous line, numerical solution with dotted line.}
    \label{fig:SingleFracture_solution}
\end{figure}

The profiles of $\|\boldsymbol{g}_T\|_2$ and $t_N$ (Figure \ref{fig:SingleFracture_solution}) are in good agreement with the analytical solution. 
The error in $t_N$ at the fault tip is a typical pathological behavior for this benchmark \cite{franceschini2020algebraically,frigo2025robust}, caused by the abrupt transition from slip to stick, which is strongly enforced due to the limited extent of the fracture. Note that the combined effect of grid non-conformity and traction-jump stabilization helps smooth out the oscillations. This error is expected to become more localized as the grid is refined.

\subsection{Test 4: Two blocks under compression and shear}
This test case was already proposed in \cite{franceschini2016novel} for conforming grids. Here, we keep the same geometrical setting and external forces. For the sake of clarity, a brief description of the model is provided below.

\begin{figure}
    \centering
    \includegraphics[width=0.7\linewidth]{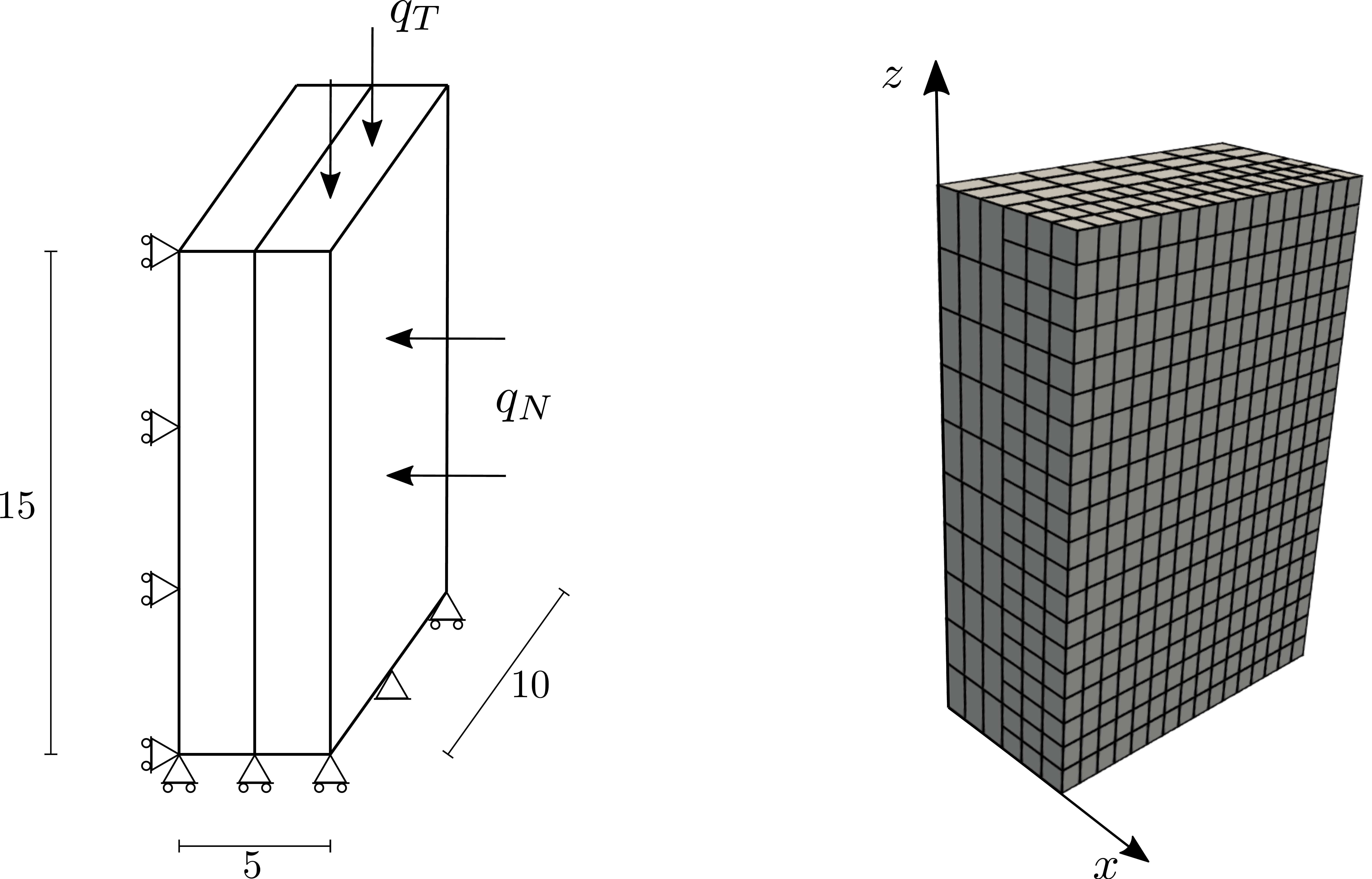}
    \caption{Test 4, two-block test case: problem setup (left) and non-conforming grids used (right).}
    \label{fig:StickSlip_setting}
\end{figure}

We consider a prismatic domain of dimensions $10 \times 5 \times 15$, composed of a linear elastic material with $E = 2000$ and $\nu = 0.25$. The prism is divided into two blocks by a vertical crack, characterized by zero cohesion and a friction angle $\phi = 30^\circ$. Boundary constraints and distributed external loads are applied as shown in Figure \ref{fig:StickSlip_setting}.

\begin{figure}
    \centering
    \includegraphics[width=0.7\linewidth]{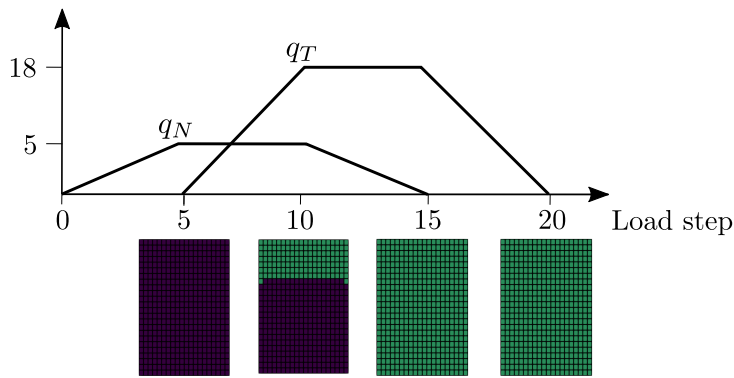}
    \caption{Test 4, two-block test case: time evolution of applied boundary traction and fracture elements in stick \filledsquare[stick]{1pt} and slip \filledsquare[slip]{1pt} state.}
    \label{fig:StickSlip_loadStep}
\end{figure}

\begin{figure}
    \centering
    \includegraphics[width=0.95\linewidth]{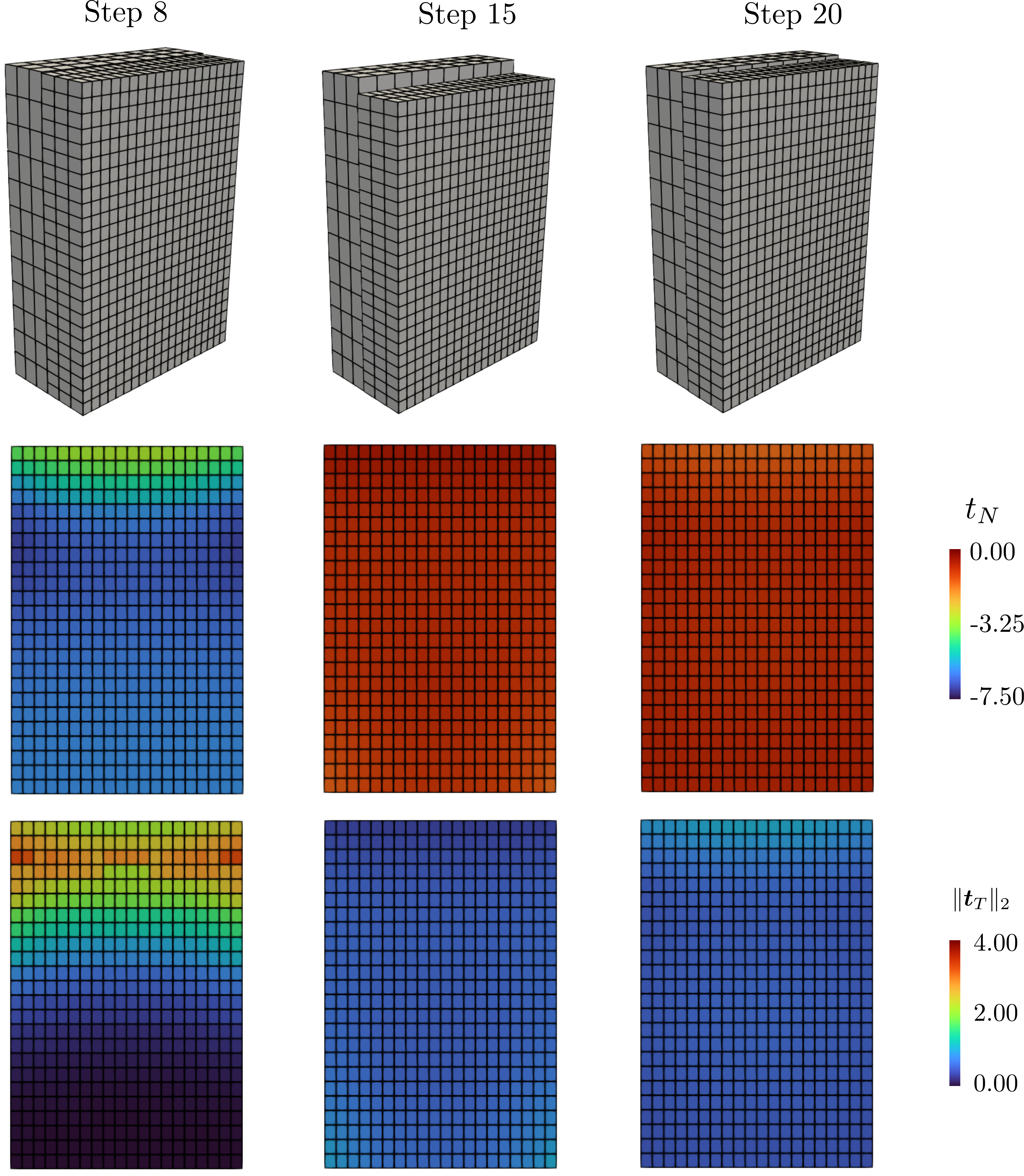}
    \caption{Test 4, two-block test case: solution at load step 8, 15 and 20. Top row: deformed configuration exaggerated by 15 times. Middle row: contour plot of the normal stress. Bottom row: contour plot of the tangential traction norm. }
    \label{fig:StickSlip_solutionEvolution}
\end{figure}

The crack is initially subjected to compressive stress $t_N = 1$. The time evolution of external normal traction $q_N$ and shear traction $q_T$ is presented in Figure~\ref{fig:StickSlip_loadStep}, along with a representation of the fracture contact states during simulation. 
On the initial loading steps, only the normal traction $q_N$ is applied and the two blocks remain in stick mode. As the shear pressure increases, the tangential traction grows until the crack starts to slip. Sliding continues as $q_N$ decreases and $q_T$ further increases. 
The simulation is carried out using a structured hexahedral grid. The mortar side is assigned to the left block, and consists of $3\times10\times12$ elements, while the non-mortar side is on the right block, made of $3\times20\times24$ elements, resulting in the mesh size ratio $r_h = 1/2$. 
%
%
Figure \ref{fig:StickSlip_solutionEvolution} reports the solution at three representative time steps. At time step 8, sliding occurs in the upper part of the interface. By time step 10, the crack has propagated, and approximately one third of the interface has entered the slip regime. At time step 15 the entire interface is sliding. The results are in full agreement with those in \cite{franceschini2016novel} for a conforming discretization.

\begin{figure}
  \centering
  \includegraphics[width=\linewidth]{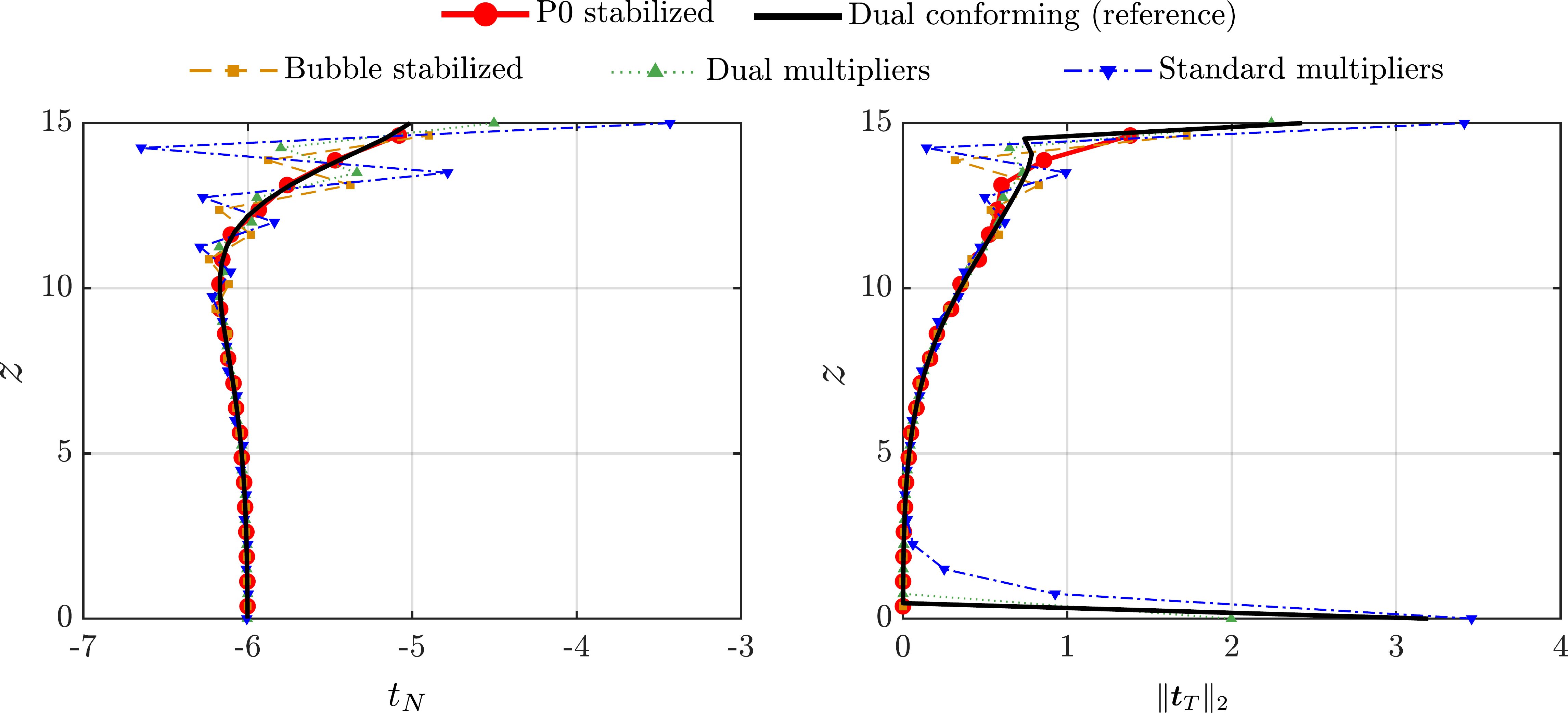}
  \caption{Test 4, two-block test case: traction distribution along the crack vertical axis at time step 6 using different multiplier formulations.}
  \label{fig:StickSlip_tractionComparison}
\end{figure}

Figure \ref{fig:StickSlip_tractionComparison} reports the traction calculated in the stick regime compared to that obtained with alternative Lagrange multiplier formulations. In particular, we consider standard and dual multipliers~\cite{popp2010dual} and the bubble-stabilized formulation introduced in~\cite{hauret2007discontinuous}. The latter employs $\mathbb{P}_0$ multipliers as in our approach, but restores the inf-sup stability condition by enriching the displacement field on the non-mortar interface with locally supported bubble functions.
As shown in~\cite{moretto2025stabilized}, when $r_h < 1$ the interface is over-constrained, thus resulting in an oscillatory traction field even if the pair of approximation spaces is inf-sup stable. In the context of Coulomb friction, numerical oscillations of the normal traction component affect the value of the limiting tangential traction, thereby compromising the robustness of the active-set strategy.
As a reference, we also include the solution obtained on a fine conforming mesh using dual multipliers. The results show that all formulations exhibit an oscillatory behavior due to the over-constraint, except the proposed traction-jump stabilized formulation that yields a smooth traction profile closely matching the reference solution.

Moreover, observe the influence on the tangential traction of the enforcement of boundary conditions at the bottom of the interface when nodal multipliers are employed.
Imposing Dirichlet boundary conditions on both mortar and non-mortar interface nodes may lead to solvability issues. To avoid this problem, we adopted a common approximate strategy \cite{puso2003mesh}, which enforces the constraint only on the mortar side, while the non-mortar side weakly inherits it. Although this approach usually works satisfactorily for the primary variable, it adversely affects the computed multipliers, i.e., the traction.
The non-local nature of standard nodal multipliers tends to propagate the boundary oscillations along the interface for several nodes before vanishing.

\subsection{Test 5: Aquifer withdrawal in a faulted domain.}
\begin{figure}
    \centering
    \includegraphics[width=\linewidth]{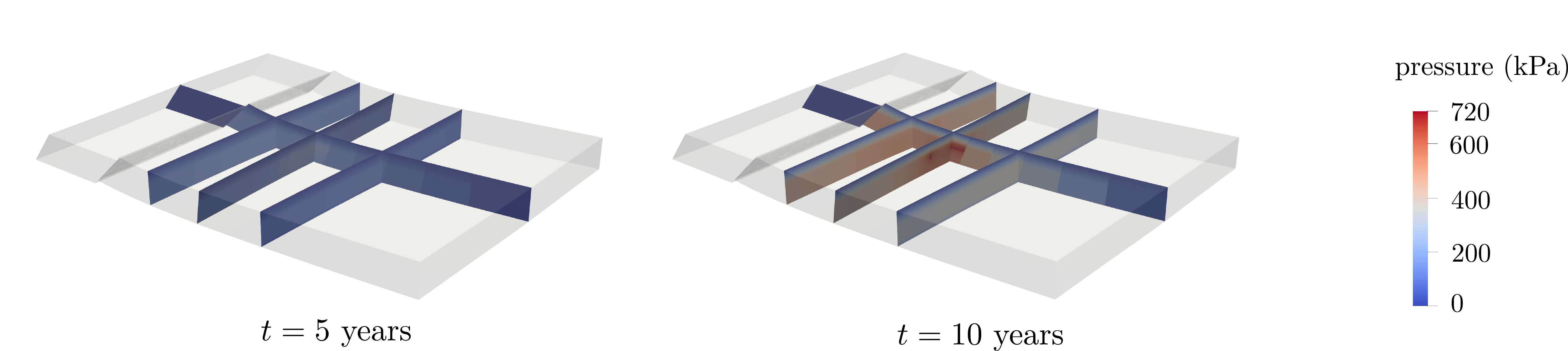}
    \caption{Test 5, aquifer withdrawal in a faulted domain: distribution of the fluid pressure after 5 and 10 years.}
    \label{fig:FaultedAquifer_pressure}
\end{figure}

The final test case is representative of a realistic application of subsurface problems in geological media. We simulate the mechanical response of a faulted geological formation due to fluid withdrawal over a period of 10 years. The problem setting, as well as the geological properties, are inspired by the test case presented in \cite{franceschini2016novel}. A passing-through discontinuity separates a rigid rock formation from a sequence of sand and clay layers. The mechanical properties selected for the model are summarized in Table \ref{tab:FaultedAquifer_materialProps}.
The fracture is tilted by $40^\circ$ with respect to the vertical axis and is characterized by a friction angle $\phi = 30^\circ$ and zero cohesion. We set an initial traction profile along the discontinuity to balance the gravitational forces.  
The fluid is extracted by two wells located in the middle of the sandy layer, at a distance of 500 m from the fault. The induced distribution of the pore pressure (Figure \ref{fig:FaultedAquifer_pressure}) is used as an external body force in the mechanical model.
The pressure change yields a contraction of the deformable sequence while the rock on the other side of the discontinuity remains practically undeformed, leading to the expected activation of the fault.

The computational grid is created with the MRST toolbox \cite{lie2019introduction} using a corner-point format, which is widely employed in industrial reservoir simulations. As is typical for such grids, the cell vertices are not aligned along the fracture interface, resulting in a non-conforming mesh that can be handled by the proposed mortar formulation. We set zero normal displacements on all external surfaces of the domain, excluding the top one. 

\begin{figure}
    \centering
    \includegraphics[width=0.8\linewidth]{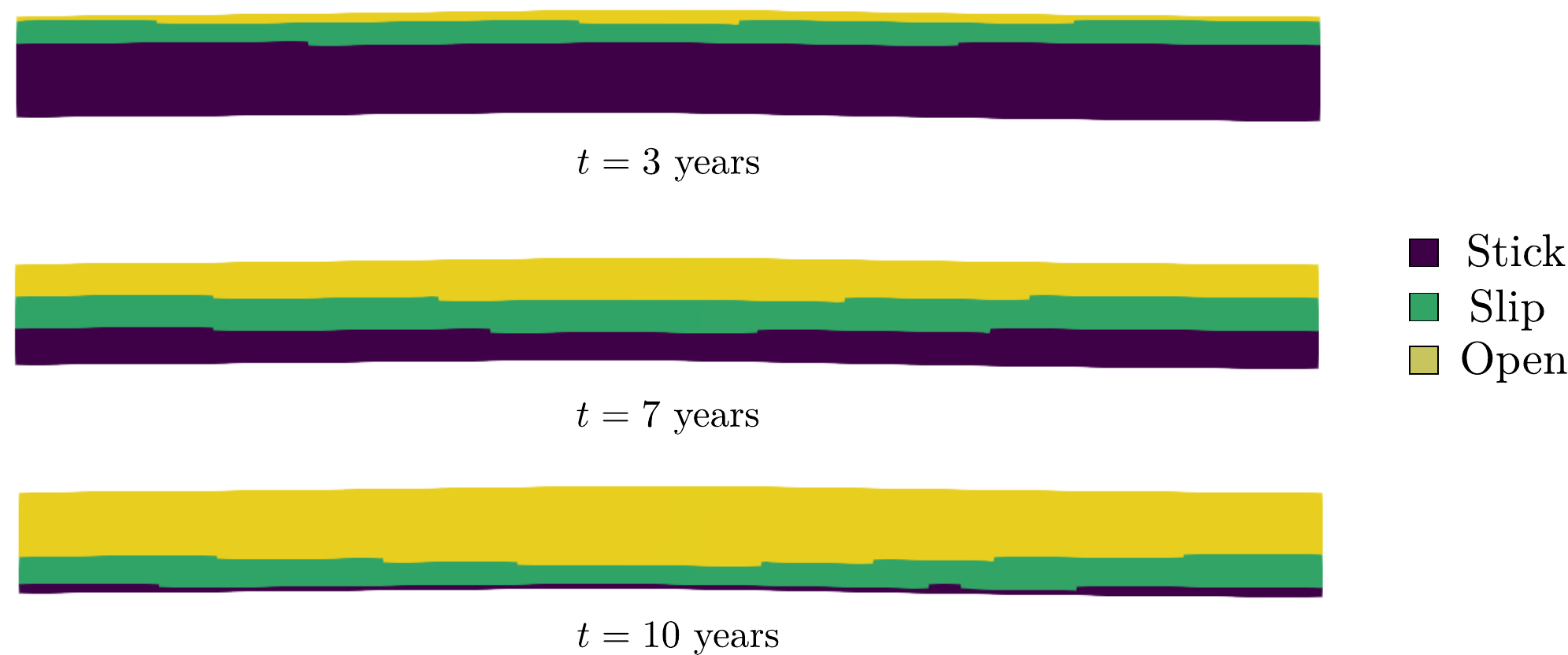}
    \caption{Test 5, aquifer withdrawal in a faulted domain: contact state of the elements located on $\Gamma_f^{(1)}$ at various time steps.}
    \label{fig:FaultedAquifer_fractureState}
\end{figure}

\begin{figure}
    \centering
    \includegraphics[width=1\linewidth]{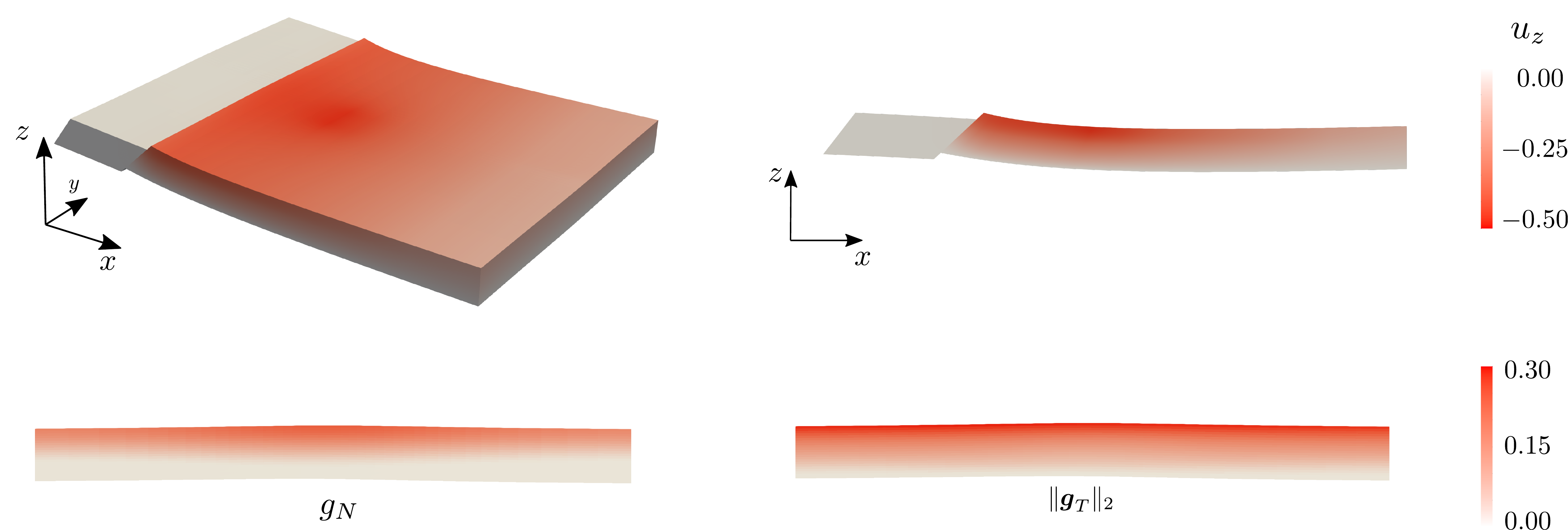}
    \caption{Test 5, aquifer withdrawal in a faulted domain: displacement results at the end of the simulation (in meters).
    Top panels: vertical displacements in the 3D domain (left) and on a slice along the x-z plane (right). Bottom panels: normal gap (left) and tangential gap (right) along the fault.}
    \label{fig:FaultedAquifer_kinematics}
\end{figure}

The resulting evolution of the fracture contact state is depicted in Figure \ref{fig:FaultedAquifer_fractureState} at three representative time steps. A sliding mechanism followed by opening initiates near the top of the domain and progressively propagates downward, until nearly the entire fracture is open at the end of the simulation. Notably, the fracture state evolves smoothly throughout the simulation.
To be consistent with the imposed displacement boundary conditions, we force the fracture elements adjacent to the bottom boundary to remain stick during the simulation.
Figure \ref{fig:FaultedAquifer_kinematics} reports the results at the end of the simulation in terms of displacement gap and global vertical displacements. The resulting gap across the discontinuity has the same order of magnitude as the maximum vertical displacement observed around the wells.

\begin{figure}
    \centering
    \includegraphics[width=\linewidth]{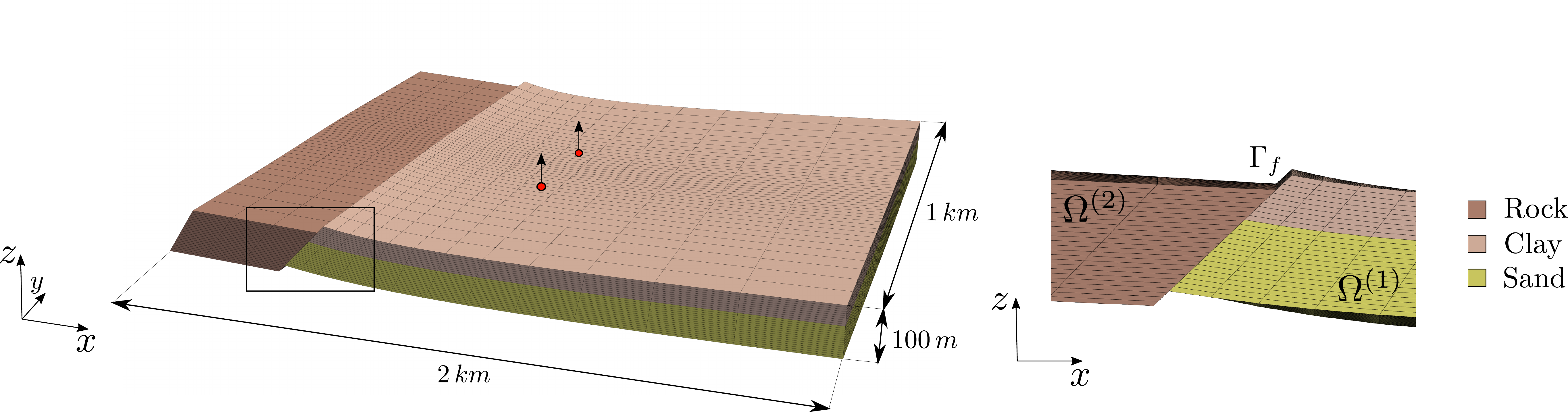}
    \caption{Test 5, aquifer withdrawal in a faulted domain: computational grid with materials definition and zoom at the fracture interface. The location of the wells is indicated by the red circles.}
    \label{fig:FaultedAquifer_setting}
\end{figure}

\begin{table}
\centering
\caption{Test 5, aquifer withdrawal in a faulted domain: material properties.}\label{tabRBF}%
\begin{tabular}{l c c}
\toprule
 &  $E$ [MPa]  & $\nu$\\
 \midrule
Rock    & 4900 & 0.25 \\ [2pt]
Clay & 30 & 0.25   \\[2pt]
Sand  & 40 & 0.25 \\[1pt]
\bottomrule
\end{tabular}
\label{tab:FaultedAquifer_materialProps}
\end{table}

\section{Conclusions}
In this work, we propose a stabilized mortar formulation for modeling contact mechanics in fractured and faulted geological media using non-conforming grids. 
Following the mortar framework, the contact constraint is enforced by a set of Lagrange multipliers defined on one side of the interface, while a Mohr-Coulomb law governs the frictional behavior of the fault. 
To resolve the variational inequality resulting from the contact constraint, we employ a standard active-set strategy.

The main contribution of this work is the use of piecewise constant ($\mathbb{P}_0$) Lagrange multipliers combined with a novel traction-jump stabilization to recover the inf-sup stability condition.
We verify numerically the inf-sup stability and develop an algorithm based on local Schur complement estimates to properly scale the inter-element traction-jump contributions, with no need for user-specified stabilization coefficients.

Numerical results demonstrate the effectiveness of the proposed stabilization in producing oscillation-free traction profiles and naturally addressing typical implementation challenges of the mortar method, such as Dirichlet boundary condition enforcement and intersecting fractures. 
Finally, we show the applicability of the proposed formulation for the simulation of contact mechanics with corner-point grids, which represent an industrial standard for many flow simulators in geological media.
The results demonstrate the flexibility and potential of the proposed approach in real-world geological settings.

Future work will focus on testing the proposed formulation in more complex geometrical scenarios and improving its robustness by leveraging more sophisticated contact algorithms compared to the standard active-set.

\vspace{0.5cm}
\noindent {\bf Acknowledgements.} M. Ferronato was supported by the RESTORE (REconstruct subsurface heterogeneities and quantify Sediment needs TO improve the REsilience of Venice saltmarshes) PRIN 2022 PNRR Project, Funded by the European Union - Next Generation EU, Mission 4, Component 1 CUP MASTER B53D2303363001. The authors are members of Gruppo Nazionale Calcolo Scientifico - Istituto Nazionale di Alta Matematica (GNCS-INdAM).









\bibliography{mybibfile}

\end{document}